\documentclass[12pt]{amsart}
\usepackage{mathpazo} 
\usepackage{avant}           
\usepackage{tensor}
\usepackage[backend=biber, style=alphabetic]{biblatex}

\AtBeginBibliography{\sloppy}
\usepackage{graphicx}
\usepackage{amssymb}

\usepackage{tikz-cd}

\newtheorem{thm}{Theorem}[section]
\newtheorem{cor}[thm]{Corollary}

\newtheorem{theorem}{Theorem}[section]

\newtheorem{lem}[thm]{Lemma}
\newtheorem{prop}[thm]{Proposition}
\theoremstyle{definition}
\newtheorem{example}[thm]{Example}
\newtheorem{defn}[thm]{Definition}
\newtheorem{conj}[thm]{Conjecture}
\theoremstyle{remark}
\newtheorem{rem}[thm]{Remark}
\numberwithin{equation}{section}

\newcommand{\supp}{\mbox{supp}}

\newcommand{\CC}{{\mathsf C}}
\newcommand{\CD}{{\mathsf D}}

\newcommand{\CF}{{\mathcal F}}

\newcommand{\A}{{A}}

\newcommand{\Z}{{\mathbb{Z}}}
\newcommand{\Q}{{\mathbb{Q}}}
\newcommand{\R}{{\mathbb{R}}}

\newcommand{\CE}{\operatorname{CE}}

\newcommand{\GL}{\operatorname{GL}}

\newcommand{\Ch}{\mathsf{Ch}}

\newcommand{\Loc}{\mathsf{Loc}}

\newcommand{\RH}{\operatorname{RH}}

\newcommand{\Hom}{\operatorname{Hom}}

\newcommand{\End}{\operatorname{End}}

\newcommand{\id}{\operatorname{Id}}

\newcommand{\rk}{\operatorname{rk}}

\newcommand{\ev}{\operatorname{ev}}

\newcommand{\pr}{{\operatorname{pr}}}

\makeatletter
\newcommand{\colim@}[2]{%
	\vtop{\m@th\ialign{##\cr
			\hfil$#1\operator@font colim$\hfil\cr
			\noalign{\nointerlineskip\kern1.5\ex@}#2\cr
			\noalign{\nointerlineskip\kern-\ex@}\cr}}%
}
\newcommand{\colim}{%
	\mathop{\mathpalette\colim@{\rightarrowfill@\textstyle}}\nmlimits@
}

\newcommand{\extp}{\@ifnextchar^\@extp{\@extp^{\,}}}
\def\@extp^#1{\mathop{\bigwedge\nolimits^{\!#1}}}

\makeatother

\begin{document}
\raggedbottom

\title{Higher Riemann--Hilbert Correspondence and Descent
for Regular Foliations}

\author{Qingyun Zeng}
\address{Department of Mathematics,
	University of Pennsylvania, Philadelphia, PA 19104}
\email{qze@math.upenn.edu}

\subjclass[2020]{Primary 53C12; Secondary 18D20, 18N60, 55U10, 58A12}
\date{}
\keywords{regular foliation, Riemann--Hilbert correspondence,
representation up to homotopy, iterated integral, infinity-local system}

\begin{abstract}
Let $\CF\subset TM$ be a regular foliation.  We construct an
$A_\infty$ integration functor from globally bounded finite-rank
leafwise cohesive modules with locally constant fiber-cohomology rank
to global plot-smooth infinity-local systems with the same regularity.
The target retains the global higher-transport objects and sheafifies
their raw Block--Smith Hom presheaves.  The construction combines the
Gugenheim--Arias Abad--Sch\"atz iterated-integral map with the
Lie-algebroid higher-holonomy formulas.  For every $M$, integration is
quasi-fully faithful.  On each star-shaped foliated box, a strictly
unital simplicial prism gives an explicit raw local inverse.  If $M$
is compact, effective descent for cohesive modules and \v{C}ech
descent for the target Hom sheaves promote the local comparison to an
$A_\infty$ quasi-equivalence on the regular full subcategories;
compact proper submersions give a distinguished class of examples.
At every fixed finite amplitude, the local comparison also induces an
equivalence of $\operatorname{Cat}_\infty$-valued
hypersheafifications.  No global raw-Hom comparison or analogous
correspondence for singular foliations or general
$L_\infty$-algebroids is asserted.
\end{abstract}

\maketitle
\enlargethispage{2pt}
\tableofcontents

\section{Introduction}

For a smooth manifold $M$, the de Rham integration map is a
quasi-isomorphism from differential forms to smooth singular cochains.
Gugenheim enhanced this comparison to an $A_{\infty}$ quasi-isomorphism
\cite{Gug77}; a convenient modern account is
\cite[Theorem 3.25]{AS12}.  Block and Smith used iterated integrals and
higher holonomy to prove a higher Riemann--Hilbert correspondence for a
compact manifold \cite[Theorem 4.1]{BS14}.

\subsection{Prior work}

For compact manifolds, Block--Smith prove the higher
Riemann--Hilbert quasi-equivalence for Block's dg-category of cohesive
modules and infinity-local systems \cite[Theorem 4.1]{BS14}.  In their proof, the
higher-holonomy formulas define an $A_\infty$ functor
\cite[Theorem 4.2]{BS14}; quasi-full faithfulness is
\cite[Proposition 4.3]{BS14}, and essential surjectivity is
\cite[Theorem 4.5]{BS14}.  Arias Abad--Sch\"atz recast this holonomy
construction through Gugenheim's $A_\infty$ de Rham map, prove the
unital object-level Lie-algebroid integration statement
\cite[Theorem 4.13]{AS12}, and construct the higher morphism
components of the unital Lie-algebroid $A_\infty$ functor
\cite[Definition 4.16 and Theorem 4.18]{AS12}.  Although the final
journal version \cite{BS14} appeared after \cite{AS12}, the latter
cites the 2009 Block--Smith preprint as its starting point.  We
therefore cite Block--Smith for the manifold superconnection and
higher-holonomy formulas, and Arias Abad--Sch\"atz for the
Lie-algebroid extension, unital target, and full higher morphism
construction.

Milnor descent for cohesive dg-categories was established by
Ben-Bassat and Block \cite{BenBasBlock13}.  Under suitable hypotheses,
the dg-category associated to a Cartesian diagram of curved dgas is
homotopy Cartesian; this yields descent for certain complex-manifold
partitions.
Block--Holstein--Wei identify twisted complexes with explicit
homotopy limits of dg-categories \cite[Theorem 4.1]{BHW17}.  Wei
subsequently proved effectivity for cohesive modules over the
Dolbeault dga \cite[Theorem 1.1]{Wei23}.  We verify that the
effectivity argument uses only the fine smooth degree-zero algebra,
finite local freeness of the positive-degree terms, perfect-complex
strictification, and Block's quasi-representability theorem.  These
inputs also hold for the leafwise de Rham dga.  This gives a
finite-cover descent theorem for leafwise cohesive modules.  The
subsequent strict Fr\'echet theory provides explicit topological
models under additional hypotheses; these hypotheses play no role in
the algebraic correspondence.

Holstein's Morita-cohomology papers give the closest dg-categorical
hyperdescent and homotopy-locally-constant-sheaf prior art.  In
\cite{Hol15Morita}, derived global sections of a constant presheaf of
dg-categories and singular cohomology are compared and identified, up
to quasi-equivalence, with perfect representations of chains on the
based loop space.  In \cite{Hol15HLC}, Morita cohomology is identified
with homotopy locally constant sheaves of perfect complexes.  These
results provide the Morita/hyperdescent and homotopy-locally-constant
sheaf context, but they do not construct the present leafwise de Rham
integration functor, the transversely smooth foliated comparison, or
the fixed-global-object sheafified-Hom target.

Chuang, Holstein, and Lazarev develop Maurer--Cartan moduli and
Riemann--Hilbert-type theorems in a homotopical algebraic framework,
including a generalization of the Block--Smith correspondence to
simplicial complexes and topological spaces \cite{CHL21}.  Their
results clarify how strong-homotopy-equivalent cochain models control
Maurer--Cartan data.  The issue addressed here is different: for a
foliation, the raw global plot-smooth Hom presheaves may
fail to be separated, while the global plot-smooth objects remain
geometrically meaningful.

In a different direction, To\"en and Vezzosi prove a
Riemann--Hilbert theorem for perfect
crystals on proper smooth complex varieties equipped with
quasi-smooth rigid derived foliations \cite{TV23}.  Their theorem
includes singular algebraic geometry and derived foliations, whereas
the present theorem concerns real smooth regular foliations,
functional smoothness on foliated simplices, and a
fixed-global-object sheafified-Hom target.  Neither result is presented
as a formal specialization of the other.

The $A_\infty$ Poincar\'e lemma of Arias Abad, Quintero V\'elez, and
V\'elez V\'asquez supplies a related local higher-homotopical
comparison \cite{AQVV19}.  Our local result additionally tracks
functional smoothness in the transverse parameter and constructs the
strictly unital raw prism used in compact globalization.  The work of
Arias Abad and V\'elez V\'asquez on principal $2$-bundles develops
higher parallel transport in a geometric $2$-categorical setting
\cite{AVV21}; it is complementary to the dg-categorical target used
here.  Recent Chern--Weil theory for infinity-local systems constructs
a categorified Chern--Weil functor \cite{AMV21}.  The compact
correspondence provides a potential leafwise source model for such
constructions; compatibility with the Chern--Weil functor lies beyond
the present scope.

Two companion manuscripts in preparation develop related
cohesive and superconnection methods in different geometries.
Zeng studies derived categories of transversely holomorphic
foliations \cite{ZengTransverseHolomorphicPrep}, while Block and Zeng
develop a categorical Penrose--Ward correspondence for self-dual
Yang--Mills theory \cite{BlockZengPenroseWardPrep}.  Neither manuscript
is used in the proofs below.

These comparisons identify the imported higher-holonomy,
Lie-algebroid integration, cohesive descent, and Morita/hyperdescent
machinery without making a priority claim for those inputs.  The new
ingredients here are the transversely smooth foliated comparison, the
fixed-global-object sheafified-Hom target, the explicit local prism
reconstruction, and the compact local-to-global argument for regular
foliations.  The global comparison between raw morphisms and their
sheafifications remains open.

\subsection{Main results and proof strategy}

A regular foliation $\CF\subset TM$ is itself a Lie algebroid, so this
extension supplies an $A_{\infty}$ integration functor from
representations up to homotopy of $\CF$ to representations up to
homotopy of its infinity-groupoid.  The purpose of this paper is to make
the leafwise model explicit, prove its de Rham comparison, construct a
local inverse by coherent prisms, and globalize that inverse by dg
descent.

Target objects carry global plot-smooth higher transport on all
foliated simplices and remain unchanged.  Only their Hom presheaves
are sheafified, so a morphism needs only local raw Block--Smith
representatives.  The comparison with the germ functor from global raw
morphisms remains open.  This fixed-object sheafified-Hom category is
the target of the compact global correspondence below; it changes
morphisms but does not add descent objects.

\begin{theorem}[Higher Riemann--Hilbert correspondence]
\label{thm:intro-main}
Let $(M,\CF)$ be a regular foliated manifold.  Let
\[
 \mathcal P_{\Omega^\bullet_{\CF}(M)}^{\mathrm{reg}}
\]
be the full dg subcategory of globally bounded finite-rank cohesive
modules whose fiberwise internal cohomology dimensions are locally
constant in every degree.  Let
\[
 \Loc_{\mathrm{shHom},\Ch_{\R}}^{\mathrm{sm,perf,reg}}
 \bigl(\Pi_\infty(\CF)\bigr)
\]
be the fixed-global-object dg category of globally bounded
finite-rank, strictly unital, normalized, plot-smooth infinity-local
systems with locally constant vertex-cohomology dimensions and
sheafified raw Block--Smith Hom complexes.  Higher integration defines
an $A_\infty$ functor
\[
 \RH_{\mathrm{shHom}}:
 \mathcal P_{\Omega^\bullet_{\CF}(M)}^{\mathrm{reg}}
 \longrightarrow
 \Loc_{\mathrm{shHom},\Ch_{\R}}^{\mathrm{sm,perf,reg}}
 \bigl(\Pi_\infty(\CF)\bigr)
\]
with the following properties.
\begin{enumerate}
 \item For arbitrary $M$, $\RH_{\mathrm{shHom}}$ is quasi-fully
 faithful.
 \item On a star-shaped foliated box, both the raw and
 sheafified-Hom integration functors are $A_\infty$
 quasi-equivalences.
 \item If $M$ is compact, $\RH_{\mathrm{shHom}}$ is an $A_\infty$
 quasi-equivalence.  This includes proper submersions over compact
 bases.
 \item After restricting both local theories to any fixed finite
 amplitude interval, applying the $A_\infty$ nerve, and
 hypersheafifying, raw local integration induces an equivalence of
 $\operatorname{Cat}_\infty$-valued hypersheaves.
\end{enumerate}
\end{theorem}

The theorem does not assert that the canonical germ map
\[
 \mathcal H^\bullet_{V,W}(M)
 \longrightarrow
 \Gamma\bigl(M,\underline{\mathcal H}^\bullet_{V,W}\bigr)
\]
is a quasi-isomorphism for global raw objects, nor does it identify the
fixed-global-object target with the global sections of an object
stackification.

The proof has four stages.  First, leafwise forms and sheafified
plot-smooth cochains are compared, with flat coefficients and then
$A_\infty$-multiplicatively.  Second, the leading arity term of the
Block--Smith generalized-cone formula is identified with coefficient
integration, giving quasi-full faithfulness.  Third, a plaque
contraction and its strictly unital prism reconstruct every target
object on a star-shaped foliated box.  Finally, effective source
descent globalizes the local inverse objects, while target \v{C}ech
full faithfulness globalizes only the equivalence between the two
objects that already exist globally.  The nuclear Fr\'echet theory of
Section~\ref{sec:topological-enhancements} is logically independent of
these four steps and furnishes strict models under additional
hypotheses.

\subsection*{Acknowledgments}

This article grew out of part of the author's Ph.D.\ thesis
\cite{ZengThesis} and forms part of a broader program applying
higher-categorical methods to foliations and related geometry.
Companion directions include characteristic classes of
$L_\infty$-resolved singular foliations
\cite{BlockZengSingularClassesPrep}, a categorical Penrose--Ward
correspondence for self-dual Yang--Mills theory
\cite{BlockZengPenroseWardPrep}, global Dolbeault
$L_\infty$-algebroid resolutions of singular holomorphic foliations
\cite{BlockWeiZengDolbeaultPrep}, and derived categories of
transversely holomorphic foliations
\cite{ZengTransverseHolomorphicPrep}.

The author is deeply grateful to Jonathan Block for his guidance,
many illuminating discussions, encouragement, and sustained support.
The author also thanks Jim Stasheff, Sylvain Lavau, and Zhiyou Wu for
helpful discussions.

\section{Conventions and foundations}

\subsection{Regular foliations and leafwise forms}

All manifolds are smooth, Hausdorff, and second countable.  Compactness is
imposed only when explicitly stated.  Complexes are cohomologically
graded and the ground field is $\R$.

\begin{defn}\label{def:foliation}
A \emph{regular foliation} on $M$ is a vector subbundle
$\CF\subset TM$ whose smooth sections are closed under the Lie bracket.
A map $f:(M,\CF_M)\to(N,\CF_N)$ is \emph{foliated} if
$df(\CF_M)\subset\CF_N$.
\end{defn}

Write
\[
 \Omega^r_{\CF}(U)=
 \Gamma\bigl(U,\extp^r\CF^\vee\bigr).
\]
The Lie algebroid structure on $\CF$ gives the leafwise differential
\begin{align*}
(d_{\CF}\omega)(X_0,\ldots,X_r)
={}&\sum_i(-1)^i X_i
  \omega(X_0,\ldots,\widehat X_i,\ldots,X_r)\\
&+\sum_{i<j}(-1)^{i+j}
  \omega([X_i,X_j],X_0,\ldots,\widehat X_i,\ldots,
  \widehat X_j,\ldots,X_r).
\end{align*}
Thus
\[
 \A=\Omega^\bullet_{\CF}(M)
\]
is a differential graded algebra.  Let
\[
 \underline{\R}_{\CF}
 =\ker\bigl(d_{\CF}:C^\infty_M\longrightarrow\Omega^1_{\CF}\bigr)
\]
be the sheaf of smooth functions that are locally constant along leaves.

\subsection{Cohesive modules}

\begin{defn}\label{def:cohesive}
A \emph{$\Z$-connection} on a bounded graded finitely generated
projective $C^\infty(M)$-module $E^\bullet$ is a total degree-one map
\[
 \mathbb E:
 E^\bullet\otimes_{C^\infty(M)}\A
 \longrightarrow
 E^\bullet\otimes_{C^\infty(M)}\A
\]
satisfying
\[
 \mathbb E(ea)=\mathbb E(e)a+(-1)^{|e|}e\,d_{\CF}a.
\]
It has components
\[
 \mathbb E^r:E^j\longrightarrow
 E^{j+1-r}\otimes_{C^\infty(M)}\Omega^r_{\CF}(M).
\]
A \emph{cohesive module} is a pair $(E^\bullet,\mathbb E)$ with
$\mathbb E^2=0$.
\end{defn}

This is Block's cohesive-module definition specialized to the
leafwise de Rham dga \cite[Definition 2.4]{Blo05}.

A \emph{quasi-cohesive module} obeys the same definition without the
bounded finite-projective condition on $E^\bullet$
\cite[Definition 3.7]{Blo05}.  Such an object need not itself lie in
$\mathcal P_{\A}$; it instead determines a dg-module over
$\mathcal P_{\A}$.  Block's second paper introduces a complementary
topological enlargement, the \emph{quasi-perfect} category: there the
underlying complex remains bounded and projective as a Fr\'echet module,
but need not be finitely generated; the connection is continuous and
relative tensor products are completed \cite[Section 2.2]{Blo06II}.  We use
``quasi-cohesive'' in the former sense and state the additional
topological hypotheses explicitly when they are needed below.
Cohesive modules form a dg-category $\mathcal P_{\A}$.  Its morphism
complex consists of $\A$-linear maps, with differential
\[
 D\phi=\mathbb F\phi-(-1)^{|\phi|}\phi\mathbb E.
\]
Shifts and mapping cones give $\mathcal P_{\A}$ its standard
pre-triangulated structure.  A cohesive module is
\emph{cohomologically regular} if every
$H^j(E^\bullet,\mathbb E^0)$ has locally constant rank, hence is the
module of sections of a vector bundle.  Denote the full dg-subcategory of
such objects by $\mathcal P_{\A}^{\mathrm{reg}}$.

\subsection{$A_{\infty}$ conventions}

We use the suspended bar convention of
\cite[Appendix A]{AS12}.  For the present paper, an $A_{\infty}$ functor
$F:\CC\to\CD$ is \emph{quasi-fully faithful} if
\[
 F_1:\CC(X,Y)\longrightarrow\CD(FX,FY)
\]
is a quasi-isomorphism for all $X,Y$.  It is
\emph{essentially surjective} if every target object is isomorphic in
$H^0(\CD)$ to an object in the image, and it is a
\emph{quasi-equivalence} if both conditions hold.

\subsection{The smooth monodromy simplicial space}

\begin{defn}\label{def:monodromy}
Regard $\Delta^n$ as a manifold with corners, and call a map from it
smooth when it extends smoothly to a neighborhood in its affine span.
For $n\geq0$, let
\[
 \Pi_{\infty}(\CF)_n
 =
 \left\{
 \sigma:\Delta^n\to M\ \middle|\
 \sigma\text{ is smooth and }d\sigma(T\Delta^n)\subset\CF
 \right\}.
\]
Faces and degeneracies are induced by precomposition.  A plot
$S\to\Pi_{\infty}(\CF)_n$, with $S$ a finite-dimensional smooth test
manifold, is declared smooth when its adjoint
$S\times\Delta^n\to M$ is smooth.  This functional diffeology is the
smooth structure used below.
\end{defn}

The underlying object is the smooth monodromy simplicial set of the
regular foliation.  It agrees with the Arias Abad--Sch\"atz simplicial
set of the Lie algebroid $\CF$: because the anchor
$\CF\hookrightarrow TM$ is injective, a Lie algebroid morphism
$T\Delta^n\to\CF$ is equivalently the differential of a foliated
simplex.

Let $C^\bullet_{\mathrm{raw}}(\CF)$ be the normalized cochain dga of
plot-smooth functions on $\Pi_{\infty}(\CF)_\bullet$.  For sheaf
arguments, let $\underline C^\bullet_{\CF}$ be the sheafification of
\[
 U\longmapsto C^\bullet_{\mathrm{raw}}(\CF|_U),
\]
and set
\[
 C^\bullet_{\mathrm{sh}}(\CF)
 =\Gamma(M,\underline C^\bullet_{\CF}).
\]
Restriction commutes with the simplicial coboundary and cup product, so
both operations sheafify.  The canonical germ map is a strict unital
dga map
\begin{equation}\label{eq:scalar-germ-map}
 \mathbf{sh}:
 C^\bullet_{\mathrm{raw}}(\CF)
 \longrightarrow C^\bullet_{\mathrm{sh}}(\CF).
\end{equation}
We make no quasi-isomorphism claim for this map.  A fixed-depth
barycentric subdivision and its prism preserve plot-smoothness, and
the required depth is locally uniform on each finite-dimensional
plot, but no globally smooth choice of depth is presently available.
The de Rham comparison below is instead proved directly for the
sheafified complex by a stalkwise Poincar\'e argument.

\begin{rem}\label{rem:mon-hol}
This paper uses \emph{monodromy}, not holonomy.  The holonomy groupoid
of a singular foliation is the Androulidakis--Skandalis construction
\cite{AS06}, a quotient that retains transverse holonomy germs.  For a
singular foliation, the foliation need not be a vector bundle and the simplicial
object in Definition~\ref{def:monodromy} is not the integration of a
universal $L_{\infty}$-algebroid.  No singular-foliation correspondence
is asserted here.  We use ``monodromy infinity-groupoid'' as conventional
terminology, but no Fr\'echet representability or smooth horn-filling
claim is needed below.
\end{rem}

\subsection{Smooth perfect infinity-local systems}

For $0\leq j\leq n$, evaluation at the $j$th vertex defines a smooth
map of diffeological spaces
\[
 \ev_j:\Pi_{\infty}(\CF)_n\longrightarrow M,
 \qquad
 \ev_j(\sigma)=\sigma(v_j).
\]

An infinity-local system on $\Pi_{\infty}(\CF)$ assigns a cochain
complex $V_x$ to each $x\in M$ and, for every $n\geq1$ and
$\sigma\in\Pi_{\infty}(\CF)_n$, a map
\[
 F_n(\sigma):
 V_{\sigma(v_n)}\longrightarrow V_{\sigma(v_0)}
\qquad\text{of degree }1-n,
\]
subject to the normalized Maurer--Cartan equation
\[
 dF+\widehat\delta F+F\cup F=0.
\]
We impose the unital normalization
\[
 F_1(s_0x)=\id_{V_x},
 \qquad
 F_n(s_i\sigma)=0\quad(n>1).
\]
See \cite[Section 2]{BS14} and \cite[Appendix A]{AS12} for the
equivalent suspended conventions.

\begin{defn}\label{def:smooth-local-system}
An infinity-local system is \emph{smooth and perfect} if its vertex
data are the fibers of a globally defined graded smooth vector bundle
\[
 V^\bullet=\bigoplus_{a\leq j\leq b}V^j\longrightarrow M
\]
of finite rank, its vertex differential is a
smooth bundle map $d_V:V^\bullet\to V^{\bullet+1}$ with $d_V^2=0$,
and the following plot condition holds.  For every finite-dimensional
plot $\varphi:S\to\Pi_{\infty}(\CF)_n$, the maps
$F_n(\varphi(s))$ assemble to a smooth degree-$1-n$ bundle map
\[
 \varphi^*\ev_n^*V\longrightarrow\varphi^*\ev_0^*V
\]
over $S$.

The object is \emph{cohomologically regular} if
$\dim H^j(V_x,d_{V,x})$ is locally constant in $x$ for every $j$; in
that case the cohomology spaces form graded vector bundles.  We also
call such an object \emph{regular}.
\end{defn}

For two such global objects $V,W$ and an open set $U\subset M$, let
$\mathcal H^k_{V,W}(U)$ be the vector space of raw normalized
plot-smooth Block--Smith morphism cochains on foliated simplices with
image in $U$.  Such a cochain has components
\[
 a_n(\sigma):V_{\sigma(v_n)}\longrightarrow W_{\sigma(v_0)}
 \quad\text{of degree }k-n.
\]
For $n=0$, this is a smooth degree-$k$ bundle map
$a_0:V^\bullet\to W^{\bullet+k}$.  For $n>0$, plot-smoothness has the
same meaning as in Definition~\ref{def:smooth-local-system}.  Put
\[
 \underline{\mathcal H}^k_{V,W}
 =
 \bigl(\mathcal H^k_{V,W}\bigr)^\#,
 \qquad
 \operatorname{Hom}^k_{\Loc_{\mathrm{shHom}}}(V,W)
 =
 \Gamma\bigl(M,\underline{\mathcal H}^k_{V,W}\bigr).
\]
The same global objects with the raw complexes
$\mathcal H^\bullet_{V,W}(M)$ form the unital dg-category
$\Loc_{\mathrm{raw},\Ch_{\R}}^{\mathrm{sm,perf}}
(\Pi_\infty(\CF))$.

\begin{prop}\label{prop:sheafified-target}
The bar-convolution differential, composition, and normalized identity
cochains make the global objects in
Definition~\ref{def:smooth-local-system} into a unital dg-category
\[
 \Loc_{\mathrm{shHom},\Ch_{\R}}^{\mathrm{sm,perf}}
 \bigl(\Pi_{\infty}(\CF)\bigr).
\]
The germ maps define a strict dg functor
\[
 \mathbf{Sh}:
 \Loc_{\mathrm{raw},\Ch_{\R}}^{\mathrm{sm,perf}}
 \bigl(\Pi_{\infty}(\CF)\bigr)
 \longrightarrow
 \Loc_{\mathrm{shHom},\Ch_{\R}}^{\mathrm{sm,perf}}
 \bigl(\Pi_{\infty}(\CF)\bigr)
\]
which is the identity on objects.
\end{prop}

\begin{proof}
The local bar differential and convolution composition commute with
restriction in $U$, so they induce maps of associated sheaves.
Associated sheafification of real vector spaces is exact: this can be
checked on stalks, which are filtered colimits.  It therefore preserves
the finite kernels defining normalization.  Since the vertex bundles
have finite global amplitudes, only finitely many simplicial components
contribute to any fixed total degree.  Explicitly, if $V$ and $W$ are
supported in $[a_V,b_V]$ and $[a_W,b_W]$, respectively, then an
arity-$n$ component of a total degree-$k$ morphism has internal degree
$k-n$ and can be nonzero only when
\[
 a_W-b_V\leq k-n\leq b_W-a_V.
\]
Thus each total degree is a finite direct sum of arity pieces and the
Hom complex vanishes below $a_W-b_V$.  No commutation of sheafification
with an infinite product is being used.

On a common refinement, two locally represented morphisms may be
composed by the raw convolution formula.  The result is independent of
the representatives because equality of germs and every bar operation
are local.  Associativity, the dg Leibniz rule, and $D^2=0$ hold on
every local raw representative and hence after sheafification.  The
normalized identity cochain is global and remains a strict unit.  The
germ maps preserve all operations, giving $\mathbf{Sh}$.
\end{proof}

\begin{rem}[Homwise universal property]
\label{rem:shhom-homwise-universal-property}
For each fixed pair of global raw objects $V,W$, the canonical map
\[
 \mathcal H^\bullet_{V,W}
 \longrightarrow
 \underline{\mathcal H}^\bullet_{V,W}
\]
has the ordinary universal property of sheafification: every morphism
from the raw Hom presheaf to a sheaf of complexes factors uniquely
through $\underline{\mathcal H}^\bullet_{V,W}$.  The differential,
composition, and units are local operations, so they descend to these
sheafifications.  In this Homwise sense the construction sheafifies
the raw morphism theory while retaining the original global object
set.  This pairwise universal property does not furnish a formal
adjunction or an initiality theorem in a $2$-category of dg or
$A_\infty$ prestacks.
\end{rem}

The full subcategory of regular objects is denoted by
\[
 \Loc_{\mathrm{shHom},\Ch_{\R}}^{\mathrm{sm,perf,reg}}
 \bigl(\Pi_{\infty}(\CF)\bigr).
\]
When the ambient foliation is clear, we abbreviate the two dg
categories to $\Loc_{\mathrm{raw}}$ and $\Loc_{\mathrm{shHom}}$.
Restriction of bundles, structure maps, and morphism germs gives dg
functors for open inclusions.  Canonically these form a coherent dg
pseudopresheaf; a fixed functorial plus-plus sheafification makes it a
strict dg-presheaf.  This statement does not assert effective descent
for the fixed global objects.
Morphisms are locally represented by plot-smooth cochains; they need
not admit one global raw representative.  The functor $\mathbf{Sh}$ is
not known to be quasi-fully faithful or conservative on object
equivalences.  In particular, no raw small-simplex comparison is
implicit in the notation.

Regularity is a restriction on objects, not a stability condition.
Even where generalized cones exist, the cone of a morphism between
regular objects can have jumping cohomology rank.  Thus neither
the regular cohesive source nor any regular target full subcategory
is asserted to be pretriangulated or stable.

\begin{prop}[Foliated pullback]\label{prop:sheafified-pullback}
Every foliated map
$f:(M,\CF_M)\to(N,\CF_N)$ induces a dg functor
\[
 f^*:
 \Loc_{\mathrm{shHom},\Ch_{\R}}^{\mathrm{sm,perf}}
 \bigl(\Pi_\infty(\CF_N)\bigr)
 \longrightarrow
 \Loc_{\mathrm{shHom},\Ch_{\R}}^{\mathrm{sm,perf}}
 \bigl(\Pi_\infty(\CF_M)\bigr).
\]
It commutes with raw pullback and the germ functors.  The same
construction gives pullback on scalar and flat-coefficient sheafified
cochains.  Pullbacks compose up to their canonical coherent
identifications, and strictly after the same fixed strictification
used for restrictions.
\end{prop}

\begin{proof}
On objects, take the graded bundle $f^*V$ and define its structure maps
by
\[
 (f^*F)_n(\sigma)=F_n(f\circ\sigma).
\]
Composition with $f$ sends plots of $\CF_M$-foliated simplices to plots
of $\CF_N$-foliated simplices, so these maps remain plot-smooth.

If $a$ is a raw morphism cochain over an open set $U\subset N$, put
\[
 (f^*a)_n(\sigma)=a_n(f\circ\sigma),
 \qquad
 \operatorname{im}(\sigma)\subset f^{-1}(U).
\]
The formula depends only on the germ of $a$ and induces a sheaf map
\[
 f^{-1}\underline{\mathcal H}^{k}_{V,W}
 \longrightarrow
 \underline{\mathcal H}^{k}_{f^*V,f^*W}
\]
on $M$.  Pulling a global section first to the inverse-image sheaf and
then applying this map defines pullback on each target Hom complex.
Local equality of representatives remains local equality after inverse
image, so this construction is independent of all representatives.

Faces, degeneracies, and endpoint evaluations commute with composition
by $f$.  Raw pullback therefore preserves the bar differential,
normalization, convolution composition, and identity cochains.  The
induced map on sheafified sections is a dg functor and commutes with
taking germs.  Scalar cochains are the same construction for the unit
object.  For coefficients, use the pulled-back flat bundle and
connection.
\end{proof}

\begin{conj}[Raw Hom comparison]\label{conj:raw-hom}
For every pair of smooth perfect global objects $V,W$, the germ map
\[
 \mathcal H^\bullet_{V,W}(M)
 \longrightarrow
 \Gamma\bigl(M,\underline{\mathcal H}^\bullet_{V,W}\bigr)
\]
is a quasi-isomorphism.  Equivalently, the identity-on-objects functor
$\mathbf{Sh}:\Loc_{\mathrm{raw}}\to\Loc_{\mathrm{shHom}}$ is a
quasi-equivalence.
\end{conj}

The scalar special case is the map
\eqref{eq:scalar-germ-map}.  The conjecture asserts only a
cohomological comparison: the raw presheaves need not be separated,
so their germ maps need not be injective on cochains.  It is also
strictly weaker than effective descent for target objects.

\begin{example}[Nonseparated raw cochains]\label{ex:nonseparated}
Take $M=\R$, $\CF=T\R$, and the cover
\[
 U=(-\infty,1),\qquad V=(-1,\infty).
\]
Choose smooth functions $f,g$ with disjoint supports concentrated near
$-2$ and $2$, respectively, and define the raw degree-one cochain
\[
 c(\sigma)=f\bigl(\sigma(v_0)\bigr)g\bigl(\sigma(v_1)\bigr).
\]
It is plot-smooth and normalized because $fg=0$.  Its restriction to
$U$ vanishes since $g|_U=0$, and its restriction to $V$ vanishes since
$f|_V=0$.  Nevertheless, $c$ is nonzero on a path from a point where
$f\neq0$ to a point where $g\neq0$.  Hence the germ map can have a
nonzero chain-level kernel; Conjecture~\ref{conj:raw-hom} says that
this defect is cohomologically harmless.
\end{example}

\section{Leafwise de Rham comparison}

\subsection{The local Poincar\'e lemma}

\begin{lem}[Leafwise Poincar\'e lemma]\label{lem:leafwise-poincare}
In a foliated chart
$U\simeq B^p\times B^q$, where the plaques are
$B^p\times\{y\}$,
\[
 \mathcal H^i(\Omega^\bullet_{\CF}|_U)
 \simeq
 \begin{cases}
 \underline{\R}_{\CF}|_U,&i=0,\\
 0,&i>0.
 \end{cases}
\]
\end{lem}

\begin{proof}
Let $h_t(x,y)=(tx,y)$ and let
$H:[0,1]\times U\to U$ be the corresponding homotopy.  For a
leafwise form $\omega$ of positive degree, set
\[
 K\omega=\int_0^1\iota_{\partial_t}H^*\omega\,dt.
\]
The standard homotopy formula gives
\[
 d_{\CF}K+Kd_{\CF}
 =\id-\pr^*\iota^*,
\]
where $\iota(y)=(0,y)$ and $\pr(x,y)=y$.  The formula depends smoothly
on the transverse parameter $y$, proving the assertion.
\end{proof}

\begin{lem}[Foliated cochain Poincar\'e lemma]\label{lem:cochain-poincare}
The augmented complex
\[
 0\longrightarrow\underline{\R}_{\CF}
 \longrightarrow\underline C^\bullet_{\CF}
\]
is locally exact.
\end{lem}

\begin{proof}
Fix a point and choose a smaller foliated box on which the
plaque-radial contraction
\[
 h_t(x,y)=(tx,y)
\]
is defined.  A foliated simplex
$\sigma:\Delta^n\to B^p\times B^q$ lies in one plaque: the differential
of its transverse component vanishes on the connected interior of
$\Delta^n$, and continuity gives the same constant transverse value on
the boundary.  For a plot of simplices this value depends smoothly on
the plot parameter, by evaluation at any fixed point of $\Delta^n$.

The homotopy $h_t$ induces a simplicial homotopy from the plaque
retraction to the identity.  Use its standard normalized prism,
obtained by summing the affine $(n+1)$-simplices associated to the
monotone maps $[n]\to[1]$.  The simplicial identities show both that
this prism preserves the degenerate chain subcomplex and that
\[
 \delta K+K\delta=\id-\pr^*\iota^*.
\]
Every prism summand is a fixed smooth affine precomposition followed
by $h$, so it preserves foliated simplices and plot-smooth families.
At $t=0$ the endpoint is the constant simplex at $(0,y)$.  Such a
simplex is degenerate in positive degree, while in degree zero the
augmentation is
\[
 c\longmapsto\bigl(y\longmapsto c(0,y)\bigr).
\]
Consequently the normalized raw cochains retract to $C^\infty(B^q)$
in degree zero and vanish in positive cohomological degrees.  This
proves exactness on the stalk at the chosen point.
No assertion that this contraction is natural under restriction to
every subopen is required.
\end{proof}

\begin{thm}[Leafwise de Rham theorem]\label{thm:leafwise-dr}
Integration over foliated simplices defines a natural quasi-isomorphism
\[
 I:\bigl(\Omega^\bullet_{\CF}(M),d_{\CF}\bigr)
 \longrightarrow
 \bigl(C^\bullet_{\mathrm{sh}}(\CF),\delta\bigr),
\qquad
 I(\omega)(\sigma)=\int_{\Delta^r}\sigma^*\omega.
\]
It is natural with respect to foliated maps.
\end{thm}

\begin{proof}
Lemmas~\ref{lem:leafwise-poincare} and
\ref{lem:cochain-poincare} show that both sheaf complexes resolve
$\underline{\R}_{\CF}$, and integration induces the identity on the
augmentation.  The sheaves of leafwise forms are fine.  Each
$\underline C^r_{\CF}$ is also fine.  Explicitly, a smooth function
$\rho$ acts in degree $r$ by
\[
 (\Theta_\rho c)(\sigma)
 =\rho\bigl(\sigma(v_0)\bigr)c(\sigma).
\]
This preserves plot-smoothness and normalization, and its
sheaf-theoretic support is contained in $\supp(\rho)$.  A locally
finite partition of unity therefore gives the required fine-sheaf
endomorphisms.  In general $\Theta_\rho$ is not a cochain map; only
termwise fineness is used.  The mapping cone of integration is a
bounded-below, sheafwise acyclic complex of fine sheaves, so its global
section complex is acyclic.  This proves the claim.
\end{proof}

\subsection{Coefficients}

Let $(L,\nabla)$ be a finite-rank vector bundle with a flat
$\CF$-connection.  A raw smooth $r$-cochain with coefficients in $L$
assigns to each foliated $r$-simplex $\sigma$ an element of
$L_{\sigma(v_0)}$, smoothly on plots.  Its differential is
\[
 (\delta_\nabla c)(\sigma)
 =
 \tau_{10}\,c(\sigma\circ\partial_0)
 +\sum_{i=1}^{r+1}(-1)^i c(\sigma\circ\partial_i),
\]
where $\tau_{10}:L_{\sigma(v_1)}\to L_{\sigma(v_0)}$ is parallel
transport along the first edge.  Sheafifying the raw plot-smooth
coefficient presheaf and taking global sections gives
$C^\bullet_{\mathrm{sh}}(\CF;L)$.  No comparison with its global raw
cochain complex is needed.  For a bounded graded flat bundle, use this
construction degree-wise and totalize.

\begin{thm}[Twisted leafwise de Rham theorem]\label{thm:twisted-dr}
Parallel transport inside each simplex to the first vertex defines a
natural quasi-isomorphism
\[
 I_L:
 \bigl(\Omega^\bullet_{\CF}(M;L),d_\nabla\bigr)
 \longrightarrow
 \bigl(C^\bullet_{\mathrm{sh}}(\CF;L),\delta_\nabla\bigr).
\]
The same conclusion holds for a globally bounded finite-rank graded
flat coefficient bundle with a parallel internal differential, after
totalization.
\end{thm}

\begin{proof}
On a star-shaped foliated chart $B^p\times B^q$, parallel transport
along the radial paths
\[
 (x,y)\longrightarrow(0,y)
\]
identifies $L$ with the pullback of
$L_0=L|_{\{0\}\times B^q}$.  Flatness makes the transport between any
two points of a plaque equal to the corresponding composite through
$(0,y)$, and smooth dependence for linear ordinary differential
equations makes this trivialization smooth in $(x,y)$.  For a graded
coefficient complex, parallelness of the internal differential makes
the trivialization a chain map.

In this trivialization, the normalized prism from
Lemma~\ref{lem:cochain-poincare} acts coefficientwise.  It commutes
with the internal differential and retracts the coefficient cochains
onto
\[
 \bigl(\Gamma(B^q,L_0),d_{L_0}\bigr).
\]
The radial covariant de Rham homotopy retracts the coefficient-form
complex onto the same transverse complex, and integration is the
identity there.  It is therefore a quasi-isomorphism on every stalk.

The coefficient-form and coefficient-cochain sheaves are termwise
fine, the latter by the operators $\Theta_\rho$ above.  Their
totalizations are bounded below.  The mapping-cone argument used in
Theorem~\ref{thm:leafwise-dr} therefore globalizes the stalkwise
comparison.
\end{proof}

\section{Parameterized iterated integrals and the
$A_{\infty}$ de Rham map}

\subsection{Vertical leafwise pullback and fiber integration}

The formulas below require differential forms only on finite-dimensional
parameter spaces.  This avoids imposing a Fr\'echet-manifold or
cotangent-bundle structure on a foliated path space.

\begin{defn}\label{def:vertical-pullback}
Let $R$ and $X$ be finite-dimensional manifolds, possibly with corners.
A smooth map
\[
 \Gamma:R\times X\longrightarrow M
\]
is \emph{vertically $\CF$-tangent} if
\[
 d\Gamma_{(r,x)}(0\oplus T_xX)
 \subseteq\CF_{\Gamma(r,x)}.
\]
For $\omega\in\Omega^k_{\CF}(M)$, its vertical pullback
$\Gamma_{\CF}^*\omega$ is the smooth $R$-family of forms on $X$
obtained by evaluating $\omega$ on the $X$-directional derivatives of
$\Gamma$.
\end{defn}

For each $r\in R$, the differential of $\Gamma_r:X\to M$ is a Lie
algebroid morphism $TX\to\CF$.  Consequently,
\begin{equation}\label{eq:vertical-pullback}
 d_X\Gamma_{\CF}^*\omega
 =\Gamma_{\CF}^*d_{\CF}\omega,
 \qquad
 \Gamma_{\CF}^*(\omega\wedge\eta)
 =\Gamma_{\CF}^*\omega\wedge\Gamma_{\CF}^*\eta.
\end{equation}
The coefficients of these forms depend smoothly on $r$.  In
particular, a plot $R\to\Pi_\infty(\CF)_n$ with adjoint
$\Sigma:R\times\Delta^n\to M$ produces a smooth $R$-family
$\Sigma_{\CF}^*\omega$ of forms on $\Delta^n$, although the
$R$-directions themselves need not be tangent to $\CF$.

Let $K$ be a compact oriented $d$-manifold with corners and let
$\pi:R\times X\times K\to R\times X$ be the projection.  Integration
over $K$ gives
\[
 \pi_!:
 C^\infty\bigl(R,\Omega^\bullet(X\times K)\bigr)
 \longrightarrow
 C^\infty\bigl(R,\Omega^{\bullet-d}(X)\bigr).
\]
With the orientation convention of \cite[Section 3.2]{AS12}, it
satisfies
\begin{align}
 \pi_!(\pi^*\beta\wedge\alpha)
 &=
 (-1)^{d|\beta|}\beta\wedge\pi_!\alpha,
 \label{eq:fiber-module}\\
 \pi_!d-(-1)^d d\pi_!
 &=
 (\pi_{\partial K})_!\iota^*.
 \label{eq:fiber-stokes}
\end{align}
Both identities reduce in local coordinates to ordinary fiberwise
integration and Stokes' theorem.  Differentiation under the integral on
compact subsets proves smooth dependence on $R$.

\subsection{Chen iterated integrals}

We follow Chen \cite{Che77}.  For the leafwise, smoothly parameterized
construction, put
\[
 \Delta_r^{\mathrm{ord}}
 =
 \{(t_1,\ldots,t_r)\mid
 1\geq t_1\geq\cdots\geq t_r\geq0\}.
\]
Suppose that
\[
 \Gamma:R\times Q\times I\longrightarrow M
\]
is smooth and tangent to $\CF$ in the $Q\times I$ directions.  Its
ordered evaluation map is
\[
 \Gamma^{(r)}:
 R\times Q\times\Delta_r^{\mathrm{ord}}\longrightarrow M^r,
 \quad
 (u,q,\mathbf t)\longmapsto
 \bigl(\Gamma(u,q,t_1),\ldots,\Gamma(u,q,t_r)\bigr).
\]
Let $sV=V[1]$ and write $[\omega]=|\omega|-1$ for reduced degree.

\begin{defn}[Chen integral]\label{def:chen}
For homogeneous leafwise forms $\omega_1,\ldots,\omega_r$, define
\[
 \operatorname{It}_{\Gamma}
 (s\omega_1|\cdots|s\omega_r)
 =
 (-1)^{\epsilon}
 \int_{\Delta_r^{\mathrm{ord}}}
 (\Gamma^{(r)})_{\CF}^*
 (\omega_1\boxtimes\cdots\boxtimes\omega_r),
\]
where
\[
 \epsilon=\sum_{i=1}^{r-1}[\omega_i](r-i).
\]
It is a smooth $R$-family of forms on $Q$ of degree
$\sum_i[\omega_i]$.
\end{defn}

Definition~\ref{def:chen} is stated for smooth families.  The only
nonsmooth family used below is Igusa's fixed piecewise-linear family;
its cellwise extension, including the required boundary regularity, is
recorded in Lemma~\ref{lem:cellwise-chen-igusa}.

If one of the $\omega_i$ is a function, the integral vanishes: the
integrand has no $dt_i$ component and hence no top component in the
ordered-simplex variables.  The construction is natural under
foliated maps and under pullback in the finite-dimensional parameter
spaces.

Let $\overline D$ be the bar differential determined by
$(\Omega^\bullet_{\CF}(M),-d_{\CF},\wedge)$.  Applying
Equation~\eqref{eq:fiber-stokes} to the faces $t_i=t_{i+1}$ and to the
two endpoint faces gives Chen's identity
\begin{align}
 d\,\operatorname{It}_{\Gamma}
 (s\omega_1|\cdots|s\omega_r)
 ={}&
 \operatorname{It}_{\Gamma}
 \bigl(\overline D(s\omega_1|\cdots|s\omega_r)\bigr)
 \notag\\
 &+\ev_1^*\omega_1\wedge
 \operatorname{It}_{\Gamma}(s\omega_2|\cdots|s\omega_r)
 \notag\\
 &-(-1)^{[\omega_1]+\cdots+[\omega_{r-1}]}
 \operatorname{It}_{\Gamma}
 (s\omega_1|\cdots|s\omega_{r-1})
 \wedge\ev_0^*\omega_r.
 \label{eq:chen-stokes}
\end{align}
The signs are precisely those of \cite[Theorem 3.10]{AS12}; the
foliation changes only the class of admissible pullbacks, not the
fiber-orientation calculation.

\begin{lem}[Reparametrization invariance]\label{lem:reparam}
Chen's iterated integral of leafwise forms is invariant under every
orientation-preserving $C^1$ diffeomorphism of $I$ fixing its
endpoints.  The same holds for piecewise-$C^1$, nondecreasing,
surjective reparametrizations.
\end{lem}

\begin{proof}
For a finite-dimensional family $\Gamma$ as above, the Chen expression
is a fiber integral over the ordered simplex.  For
$u\in\operatorname{Diff}^+(I)$, the map
\[
 (t_1,\ldots,t_n)\longmapsto
 (u(t_1),\ldots,u(t_n))
\]
is an orientation-preserving diffeomorphism of the ordered simplex
$1\geq t_1\geq\cdots\geq t_n\geq0$.  Change of variables for fiber
integration therefore gives equality as smooth $R$-families of
differential forms on $Q$.
The argument is plaque-wise and smooth in the transverse parameter.
For a non-strict monotone map the same identity follows by applying the
one-dimensional substitution formula successively to the ordered
integral; flat intervals contribute zero.  Applying that argument on
each interval of a finite partition proves the piecewise-$C^1$ case.
\end{proof}

\begin{rem}
The last statement is not asserted for arbitrary monotone continuous
maps.  Singular reparametrizations need not satisfy the substitution
formula used in the proof.
\end{rem}

\subsection{Igusa's cube-to-simplex families}

For $k\geq1$, Igusa constructs a piecewise-smooth family
\[
 \Theta_{(k)}:
 I^{k-1}\longrightarrow P(\Delta^k;v_k,v_0)
\]
of paths from the last vertex to the first
\cite[Proposition 4.6]{Igu09}.  Its adjoint is
\[
 \Theta_k:I^{k-1}\times I\longrightarrow\Delta^k.
\]
The family is obtained by sending the cube monotonically to the ordered
simplex and joining its successive coordinate vertices.  Two boundary
properties contain the geometry needed below.

\begin{lem}[Igusa fundamental-chain normalization]
\label{lem:igusa-degree}
With the orientation convention used here, the adjoint
$\Theta_k:I^k\to\Delta^k$ has degree $(-1)^k$.  Thus, for every
$k$-form $\alpha$ on $\Delta^k$,
\[
 \int_{I^k}\Theta_k^*\alpha
 =
 (-1)^k\int_{\Delta^k}\alpha.
\]
The same identity holds for a flat coefficient bundle after parallel
trivialization over $\Delta^k$.
\end{lem}

\begin{proof}
The degree calculation for Igusa's explicit piecewise-linear map is
\cite[Remark 3.17(2)]{AS12}; it is also implicit in
\cite[Proposition 4.6]{Igu09}.  A flat bundle on the contractible
simplex has a parallel trivialization.  In that frame, the coefficient
identity is the displayed scalar identity applied componentwise.
\end{proof}

\begin{prop}[Igusa face identities]\label{prop:igusa-faces}
For $1\leq i\leq k-1$:
\begin{enumerate}
 \item on the negative cubical face $w_i=0$, $\Theta_{(k)}$ is
 the path family $\Theta_{(k-1)}$ for the $i$-th simplicial face,
 up to an endpoint-preserving piecewise-smooth reparametrization;
 \item on the positive cubical face $w_i=1$, $\Theta_{(k)}$ is the
 concatenation of the path families for the corresponding front and
 back faces of $\Delta^k$.
\end{enumerate}
\end{prop}

\begin{proof}
These are \cite[Lemmas 4.7 and 4.8]{Igu09}.  The first identity may be
used inside a Chen integral by Lemma~\ref{lem:reparam}.  For the second,
decomposing the ordered configuration simplex according to how many
times lie in each path segment gives the factorization
\[
 \operatorname{It}_{\alpha\star\beta}
 (s\omega_1|\cdots|s\omega_r)
 =
 \sum_{\ell=0}^r
 \operatorname{It}_{\alpha}
 (s\omega_1|\cdots|s\omega_\ell)
 \wedge
 \operatorname{It}_{\beta}
 (s\omega_{\ell+1}|\cdots|s\omega_r),
\]
with the suspended signs fixed by the convention of
\cite[Lemma 3.19]{AS12}.
\end{proof}

\begin{lem}[Cellwise Chen integration for the Igusa family]
\label{lem:cellwise-chen-igusa}
Fix $k,r\geq1$.  There is a finite polyhedral subdivision
$\mathcal K_{k,r}$ of
\[
 I^{k-1}\times\Delta_r^{\mathrm{ord}}
\]
such that every ordered evaluation map
\[
 (u,t_1,\ldots,t_r)\longmapsto\Theta_k(u,t_i),
 \qquad 1\leq i\leq r,
\]
is smooth on each closed cell.  The subdivision is independent of
external plot parameters.  There is also a finite polyhedral
subdivision $\mathcal Q_{k,r}$ of $I^{k-1}$ such that, for every
closed cell $D\in\mathcal Q_{k,r}$, the cells induced on
$D\times\Delta_r^{\mathrm{ord}}$ are compact manifolds with corners
whose projections to $D$ are smooth on closed cells.

Consequently, if
$\Sigma:R\times\Delta^k\to M$ is a smooth finite-dimensional plot of
foliated simplices, the Chen pullbacks along
$\Sigma\circ(\id_R\times\Theta_k)$ are smooth up to the faces of every
closed cell and agree on common faces.  Fiber integration over the
ordered-simplex variables, performed cell by cell over each
$D\times\Delta_r^{\mathrm{ord}}$, defines a smooth $R$-family of
smooth forms on $D$, hence a piecewise-smooth form on $I^{k-1}$.
Chen--Stokes and cubical Stokes hold cellwise, and all internal
codimension-one contributions cancel in pairs.
\end{lem}

\begin{proof}
Igusa's adjoint $\Theta_k$ is a finite piecewise-linear map
\cite[Proposition 4.6]{Igu09}.  For fixed $r$, take a common
polyhedral refinement subordinate to the finitely many ordered
evaluation maps; this gives $\mathcal K_{k,r}$.  Triangulating the
finite piecewise-linear projection
\[
 \mathcal K_{k,r}\longrightarrow I^{k-1}
\]
simultaneously in source and target gives $\mathcal Q_{k,r}$ and the
stated projection-compatible refinement.  On each closed source cell,
the ordered evaluations are smooth, so their composites with
$\Sigma$ are smooth up to every face and depend smoothly on $R$.
Restrictions from adjacent cells agree on their common face.
Fiber integration on compact cells therefore commutes with
differentiation in $R$, and Stokes' theorem applies cellwise.  Every
internal codimension-one face occurs twice with opposite incidence
orientations and equal restricted pullback, so the two contributions
cancel.  Only the exterior ordered-simplex and cubical faces remain.
\end{proof}

This argument uses the fixed polyhedral structure of Igusa's map.
Path reparametrization alone would not justify smoothing arbitrary
nonsmooth dependence on the other cube variables.

For a foliated $k$-simplex $\sigma$, composition with $\Theta_{(k)}$
gives a $(k-1)$-parameter family of foliated paths.  For every Chen
expression $\alpha$, set
\[
 \mathsf S(\alpha)(\sigma)
 =
 \int_{I^{k-1}}
 \Theta_{(k)}^*(P\sigma)^*\alpha.
\]
The negative cubical faces give the interior simplicial faces
$\partial_i$, $1\leq i\leq k-1$.  On a positive cubical face
$w_i=1$, Chen concatenation gives a sum over all deconcatenations of
the word.  The terms with both sides nonempty are the cup-product
terms.  If one side is empty, that factor is the degree-zero identity.
Such a term survives only when the cube factor attached to the empty
side is zero-dimensional.  Thus the empty initial split survives only
for $i=1$ and gives $\partial_0$, while the empty terminal split
survives only for $i=k-1$ and gives $\partial_k$.  For $k\geq2$, these
extreme terms and the negative cubical faces form the full simplicial
coboundary.  Thus
\begin{align}
 \mathsf S(d\alpha)
 =
 b_1'\mathsf S(\alpha)
 +
 \sum_{\substack{\alpha=\alpha'\mid\alpha''\\
                  \alpha',\alpha''\neq\varnothing}}
 b_2'\bigl(\mathsf S(\alpha')\otimes
 \mathsf S(\alpha'')\bigr)
 \label{eq:igusa-stokes}
\end{align}
for $\alpha$ in the image of the Chen map.  This is the
Chen--Igusa boundary identity of \cite[Proposition 3.22]{AS12}.
For $k=1$ the parameter cube is a point and has no cubical boundary,
so the two outer faces are supplied by the endpoint terms in
Chen--Stokes together with the separately defined unary component.
For $k=0$ the construction is unary evaluation on vertices.

\subsection{The leafwise $A_{\infty}$ de Rham theorem}

Define raw maps
\[
 \widetilde\phi_n:
 \bigl(s\Omega^\bullet_{\CF}(M)\bigr)^{\otimes n}
 \longrightarrow
 sC^\bullet_{\mathrm{raw}}(\CF)
\]
as follows.  The unary component is
\[
 \widetilde\phi_1(s\omega)(\sigma)
 =(-1)^k\int_{\Delta^k}\sigma^*\omega
 \qquad(\omega\in\Omega^k_{\CF}(M)).
\]
For $n>1$, let
\[
 \widetilde\phi_n=\mathsf S\circ\operatorname{It}.
\]
Put $\phi_n=s\mathbf{sh}\circ\widetilde\phi_n$.
The separate unary definition is needed in degree zero, since a Chen
integral containing a function vanishes whereas
$\widetilde\phi_1(f)(x)=f(x)$.

\begin{thm}[$A_{\infty}$ leafwise de Rham theorem]
\label{thm:ainfty-dr}
The maps $\{\phi_n\}_{n\geq1}$ define a normalized unital
$A_{\infty}$ quasi-isomorphism
\[
 \phi:
 \bigl(\Omega^\bullet_{\CF}(M),-d_{\CF},\wedge\bigr)
 \longrightarrow
 \bigl(C^\bullet_{\mathrm{sh}}(\CF),\delta,\cup\bigr).
\]
Naturality holds for foliated maps.
\end{thm}

\begin{proof}
Equation~\eqref{eq:chen-stokes} converts the source bar differential
into the differential of the Chen expression plus its two endpoint
terms.  Equation~\eqref{eq:igusa-stokes} converts that differential
into the target coboundary and the nonempty cup-product
decompositions.  The endpoint terms supply the products in which the
separately defined unary component occurs.  If a nonempty Chen word
contains a zero-form input, it vanishes.  Mixed identities involving
functions are instead obtained by applying Chen--Stokes to the
corresponding source-bar terms in which the function is differentiated;
for example, the identity for $(f,\omega)$ uses the endpoint identity
for $\widetilde\phi_2(df,\omega)$.  This is the suspended
$A_\infty$ identity of \cite[Theorem 3.25]{AS12}.  The remaining
zero-dimensional case, two function inputs, contains no nonempty Chen
word and reduces directly to
\[
 \widetilde\phi_1(fg)
 =
 \widetilde\phi_1(f)\cup\widetilde\phi_1(g)
\]
on vertices.

It remains to check the foliation-specific smoothness.  Let
$R\to\Pi_\infty(\CF)_k$ be a plot with adjoint
$\Sigma:R\times\Delta^k\to M$.  Definition~\ref{def:vertical-pullback}
makes each $\Sigma_{\CF}^*\omega_i$ a smooth $R$-family of forms on
$\Delta^k$.  Lemma~\ref{lem:cellwise-chen-igusa} makes the Igusa
integrand piecewise smooth on a subdivision independent of $R$.
Wedge products and cellwise integration over the compact ordered
simplices and cubes preserve smooth dependence on $R$.  Hence every
$\widetilde\phi_n$ is plot-smooth.  Degenerate simplices give the
normalization by \cite[Proposition 3.26]{AS12}, and all operations
commute with pullback by a foliated map.

The formulas first take values in the raw cochain dga; compose with
the strict germ map~\eqref{eq:scalar-germ-map}.  The $A_\infty$
identities and units survive because this map preserves the
differential and product.  Finally, $\phi_1$ is the signed version of
Theorem~\ref{thm:leafwise-dr}, so the $A_\infty$ morphism is a
quasi-isomorphism.
\end{proof}

The sign in Theorem~\ref{thm:ainfty-dr} is appropriate for the scalar
dga with differential $-d_{\CF}$.  Cohesive Hom complexes use the
opposite convention.  For a graded vector space $V$, let
\[
 J:
 \bigl(\Omega^\bullet_{\CF}(\End V),d_{\CF},\wedge\bigr)
 \longrightarrow
 \bigl(\Omega^\bullet_{\CF}(\End V),-d_{\CF},\wedge\bigr),
 \qquad
 J(\eta)=(-1)^{|\eta|}\eta,
\]
where $|\eta|$ is total degree, including the endomorphism degree.
Then $J$ is a strict dga isomorphism: it is multiplicative and
$(-d_{\CF})J=Jd_{\CF}$.  If $\psi^-$ denotes the local
endomorphism-valued version of the $A_\infty$ map in
Theorem~\ref{thm:ainfty-dr}, put
\begin{equation}\label{eq:categorical-de-rham-conjugation}
 \widetilde\psi=\psi^-\circ J.
\end{equation}
After Maurer--Cartan twisting, this map naturally has the signed
target convention in which the vertex differential associated to
$(V,d_V)$ is $-d_V$.  To return to the standard Block--Smith target,
let
\[
 S_V|_{V^i}=(-1)^i\id_{V^i}.
\]
Conjugation by these grading involutions defines a strict dg
equivalence $\mathfrak G$ from the signed convention to the standard
one.  On an arity-$k$ morphism component $a_k$ of internal degree $q$
and on an object structure component $\widetilde F_k$, respectively,
\begin{equation}\label{eq:target-grading-conjugation}
 \mathfrak G(a_k)
 =
 S_W\circ a_k\circ S_V^{-1}
 =
 (-1)^q a_k,
 \qquad
 \mathfrak G(\widetilde F_k)
 =
 (-1)^{1-k}\widetilde F_k.
\end{equation}
The formulas are local and therefore define the same strict
equivalence on the raw and sheafified-Hom targets.

This source conjugation followed by target grading conjugation is the
change of convention in \cite[Definition 4.7]{AS12}.  A
total-degree-one Maurer--Cartan form satisfies
$J(\omega)=-\omega$.  A morphism component of form degree $p$ and
internal degree $q$ first satisfies
\begin{equation}\label{eq:categorical-unary-sign}
 \begin{aligned}
 \widetilde\RH_1(\eta^{(p,q)})
 &=
 (-1)^{p+q}(-1)^p I(\eta^{(p,q)})\\
 &=
 (-1)^q I(\eta^{(p,q)}),\\
 \RH_1(\eta^{(p,q)})
 &=
 \mathfrak G\widetilde\RH_1(\eta^{(p,q)})
 =
 I(\eta^{(p,q)}).
 \end{aligned}
\end{equation}
Thus the final arity-zero map is the identity on the ordinary vertex
Hom complex, and for the unit cohesive object the edge formula is the
unsigned identity
$I(df)(\sigma)=f(\sigma(v_1))-f(\sigma(v_0))$.

\section{Integration of cohesive modules}

\subsection{Endomorphism-valued iterated integrals}

Let $(E^\bullet,\mathbb E)$ be a cohesive module and let
$\sigma:\Delta^k\to M$ be a foliated simplex.  Choose a graded
trivialization of $\sigma^*E$.  In this trivialization the pulled-back
$\Z$-connection has the form
\[
 d+\omega,\qquad
 \omega\in
 \Omega^\bullet(\Delta^k;\End(E))^1,
\]
where the superscript denotes total degree.  The trivialization is
essential: endomorphisms of $E_{\gamma(t_i)}$ at different points of a
path cannot be composed without first identifying the fibers.

Matrix multiplication, wedge product, and
Definition~\ref{def:chen} give endomorphism-valued iterated integrals.
Block--Smith write the same connection as $d-\alpha$, where
$\alpha=-\omega$ \cite[Section 3.4]{BS14}.  Their total holonomy is
therefore
\begin{equation}\label{eq:dyson}
 \Psi_\omega
 =
 \id+
 \sum_{r\geq1}
 \operatorname{It}
 (\underbrace{s\alpha|\cdots|s\alpha}_{r\ {\rm factors}})
 =
 \id+
 \sum_{r\geq1}(-1)^r
 \operatorname{It}
 (\underbrace{s\omega|\cdots|s\omega}_{r\ {\rm factors}}).
\end{equation}
This is the intermediate object-level specialization of
Equation~\eqref{eq:categorical-de-rham-conjugation}: since $\omega$
has total degree one, $J(\omega)=-\omega$.  The target grading
conjugation is applied after the Igusa integral in
Equation~\eqref{eq:standard-higher-holonomy}.  For a line connection
$d+a(t)\,dt$ along a parameterized interval,
Equation~\eqref{eq:dyson} is $\exp(-\int_0^1a(t)\,dt)$.

This is generally an infinite series.  On a compact parameter set, its
$r$-th term and every fixed parameter derivative have estimates of the
form
\[
 \left\|\operatorname{It}_r\right\|
 \leq \frac{C^r(1+r)^N}{r!}.
\]
Consequently, the series and all derivatives converge uniformly on
compact subsets of every finite-dimensional plot.  It therefore
defines a smooth form-valued homomorphism from the fiber at the endpoint
to the fiber at the starting point.

\begin{prop}[Curvature insertion identity]\label{prop:curvature}
Let
\[
 \kappa=(d+\omega)^2=d\omega+\omega\wedge\omega
\]
be the curvature of the pulled-back $\Z$-connection.  Then
\begin{equation}\label{eq:curvature-insertion}
 d\Psi_\omega
 =
 -\ev_1^*\omega\circ\Psi_\omega
 +\Psi_\omega\circ\ev_0^*\omega
 +\mathcal R_\omega,
\end{equation}
where
\[
 \mathcal R_\omega
 =
 \sum_{a,b\geq0}(-1)^{a+b}
 \operatorname{It}
 \bigl(
  (s\omega)^{\otimes a}\mid
  s\kappa\mid
  (s\omega)^{\otimes b}
 \bigr),
\]
with the remaining Koszul signs fixed by the suspended
endomorphism-valued convention.  In particular, if $\kappa=0$, then
\[
 d\Psi_\omega
 =
 -\omega|_{v_0}\circ\Psi_\omega
 +\Psi_\omega\circ\omega|_{v_k}
\]
on the path family with fixed endpoints, with the Koszul signs
understood in the suspended convention.
\end{prop}

\begin{proof}
Apply the endomorphism-valued version of
Equation~\eqref{eq:chen-stokes} to each summand of
Equation~\eqref{eq:dyson}.  A differential insertion in a word of
length $a+b+1$ and a product insertion in a word of length $a+b+2$
carry opposite Dyson signs; together they give
$d\omega+\omega^2=\kappa$ with the coefficient $(-1)^{a+b}$.
The same alternation reverses the two endpoint signs.  The estimates
above justify termwise differentiation and rearrangement.  This is
the calculation of \cite[Section 3.4]{BS14} after the change of
variable $\alpha=-\omega$; the equivalent endomorphism-valued
$A_\infty$ formula is used in \cite[Sections 4.1--4.3]{AS12}.
\end{proof}

Changing the trivialization by a degree-zero gauge
$g:\Delta^k\to\GL(E)$ replaces
\[
 \omega\quad\text{by}\quad
 \omega^g=g^{-1}\omega g+g^{-1}dg.
\]
Gauge naturality of the $A_\infty$ de Rham map
\cite[Section 4.1 and Proposition 4.4]{AS12} gives
\begin{equation}\label{eq:gauge-holonomy}
 \widetilde{\operatorname{Hol}}_k(\omega^g)
 =
 g(v_0)^{-1}\widetilde{\operatorname{Hol}}_k(\omega)g(v_k).
\end{equation}
Degree-zero gauge maps commute with the grading involutions, so the
same identity holds after applying $\mathfrak G$.  Thus the endpoint
maps defined below are independent of the trivialization.

\subsection{Higher holonomy}

For $k\geq1$, regard $\sigma$ as the Igusa $(k-1)$-parameter family of
paths from $v_k$ to $v_0$ and define
\begin{equation}\label{eq:higher-holonomy}
 \widetilde{\operatorname{Hol}}_k(\sigma;\mathbb E)
 =
 \int_{I^{k-1}}
 \Theta_{(k)}^*(P\sigma)^*\Psi_\omega:
 E_{\sigma(v_k)}\longrightarrow E_{\sigma(v_0)}.
\end{equation}
It has internal degree $1-k$.  At a vertex $x$, put
\[
 \widetilde{\operatorname{Hol}}_0(x;\mathbb E)
 =
 (E_x,-\mathbb E_x^0).
\]
These are the structure components in the signed target convention.
The final Block--Smith components are
\begin{equation}\label{eq:standard-higher-holonomy}
 \operatorname{Hol}_k(\sigma;\mathbb E)
 =
 \mathfrak G\widetilde{\operatorname{Hol}}_k(\sigma;\mathbb E)
 =
 (-1)^{1-k}
 \widetilde{\operatorname{Hol}}_k(\sigma;\mathbb E),
 \qquad k\geq1,
\end{equation}
and at a vertex
\[
 \operatorname{Hol}_0(x;\mathbb E)
 =
 (E_x,\mathbb E_x^0).
\]
Equation~\eqref{eq:gauge-holonomy} makes these definitions intrinsic.
For $k=1$, the degree-zero component is ordinary parallel transport.
For $k>1$, the higher components record the homotopies that make that
transport coherent.

\begin{prop}\label{prop:holonomy-mc}
If $\mathbb E^2=0$, the maps
$\operatorname{Hol}_k(-;\mathbb E)$ satisfy the normalized
Maurer--Cartan equations of an infinity-local system.
\end{prop}

\begin{proof}
Flatness removes the curvature-insertion terms in
Proposition~\ref{prop:curvature}.  Integrating the remaining identity
over $I^{k-1}$ and using Proposition~\ref{prop:igusa-faces} turns its
negative faces into
$\widehat\delta\widetilde{\operatorname{Hol}}$ and its positive faces
into
$\widetilde{\operatorname{Hol}}\cup
\widetilde{\operatorname{Hol}}$ in the signed convention.  Applying
the strict dg equivalence $\mathfrak G$ gives
\[
 d\operatorname{Hol}
 +\widehat\delta\operatorname{Hol}
 +\operatorname{Hol}\cup\operatorname{Hol}=0.
\]
The degeneracy calculation is
the unital refinement in \cite[Proposition 3.26]{AS12}; the complete
object-level higher-holonomy calculation is
\cite[Section 4.1]{BS14}, and its unital Lie-algebroid formulation is
\cite[Theorem 4.13]{AS12}.  Notice that only the forward implication
is needed; no converse reconstructing curvature from higher holonomy
is asserted.
\end{proof}

\subsection{The integration functor}

A cohesive module over $\A=\Omega^\bullet_{\CF}(M)$ is equivalently a
representation up to homotopy of the Lie algebroid $\CF$
\cite[Definition 3.1]{AC11}.  Under this identification, higher holonomy
along a foliated simplex defines the following functor.

\begin{thm}[Higher integration functor]
\label{thm:integration}
There is a natural $A_{\infty}$ functor
\[
 \RH:
 \mathcal P_{\A}
 \longrightarrow
 \Loc_{\mathrm{shHom},\Ch_{\R}}^{\mathrm{sm,perf}}
 \bigl(\Pi_{\infty}(\CF)\bigr).
\]
On objects, $\RH(E,\mathbb E)$ assigns to a foliated simplex the higher
holonomy of $\mathbb E$.  The unary morphism component is the
Block--Smith cone holonomy \cite[Sections 3.5 and 4.1]{BS14}; all
higher morphism components and their extension to representations up
to homotopy of the Lie algebroid $\CF$, with a unital target, are
\cite[Definition 4.16 and Theorem 4.18]{AS12}.  The morphism
components first use
$\widetilde\psi=\psi^-\circ J$ from
Equation~\eqref{eq:categorical-de-rham-conjugation} and then apply the
target grading equivalence $\mathfrak G$ from
Equation~\eqref{eq:target-grading-conjugation}.
\end{thm}

\begin{proof}
Block--Smith construct the object holonomy and the unary cone map in
\cite[Sections 3.4--4.1]{BS14}.  Arias Abad--Sch\"atz identify these
formulas with the endomorphism-valued $A_\infty$ de Rham map and
construct all higher morphism components for representations up to
homotopy of every Lie algebroid
\cite[Definition 4.16 and Theorem 4.18]{AS12}.  Pulling a cohesive
$\CF$-module back along a foliated simplex gives precisely the flat
$\Z$-connection on a simplex to which those formulas apply.  Applying
their construction to $\CF$ and using Definition~\ref{def:monodromy}
gives the displayed simplicial target.

For a plot
\[
 S\longrightarrow\Pi_\infty(\CF)_n,
\]
work locally on $S$ and trivialize the pulled-back graded bundle over
$S\times\Delta^n$.  The object and unary holonomy formulas of
\cite[Sections 3.4--4.1 and Theorem 4.2]{BS14} and the higher
components of \cite[Theorem 4.18]{AS12} are uniformly convergent
iterated-integral series of the type considered above.
The estimates after \eqref{eq:dyson} control every parameter
derivative.  Wedge product and fiber integration preserve smooth
parameter dependence, so the values vary smoothly on $S$.  The source
bundle has finite amplitude and rank, while the Arias
Abad--Sch\"atz extension supplies the unital normalization.  Hence the
resulting object and morphism cochains lie in the stated smooth perfect
raw category.  Composing the morphism components with the strict germ
functor $\mathbf{Sh}$ of Proposition~\ref{prop:sheafified-target}
preserves the $A_\infty$ identities and units and gives the stated
target.
\end{proof}

\subsection{Generalized cones and the higher functor components}

The higher maps of $\RH$ are not added formally after defining its
object map; they are the twisted components of the endomorphism-valued
$A_\infty$ de Rham map.  Let
\[
 \phi_i:(E_{i-1},\mathbb E_{i-1})
 \longrightarrow(E_i,\mathbb E_i),
 \qquad 1\leq i\leq n,
\]
be composable homogeneous morphisms.  Following
\cite[Definition 4.16]{AS12}, pull the objects and morphisms back to a
simplex and choose a graded trivialization.  Write $\omega_i$ for the
diagonal Maurer--Cartan forms corresponding to $d+\omega_i$ and
$\eta_i$ for the off-diagonal forms corresponding to $\phi_i$ under
\cite[Proposition A.5]{AS12}.  On
$V=\bigoplus_{i=0}^n V_i$, put
\[
 \Omega=
 \begin{pmatrix}
  \omega_0&0&\cdots&0\\
  \eta_1&\omega_1&\ddots&\vdots\\
  0&\eta_2&\ddots&0\\
  \vdots&\ddots&\ddots&\omega_n
 \end{pmatrix}.
\]
For the categorical map
\eqref{eq:categorical-de-rham-conjugation}, set
\[
 \alpha_i=J(\omega_i)=-\omega_i,
 \qquad
 \eta_i^J=J(\eta_i)=(-1)^{|\phi_i|}\eta_i,
\]
before applying the shifts dictated by the degrees of the $\phi_i$.
The shifted-block Koszul signs are those of the suspended generalized
cone.  Once the $J$-signs and the shifted-block signs have been
absorbed into $\eta_i^J$, a length-$r$ word in the $[n,0]$ corner has
the shorthand coefficient $(-1)^{r-n}$, one minus sign for each
diagonal connection input.  If the untransformed bare $\eta_i$ are
used instead, the additional factor
$(-1)^{\sum_i|\phi_i|}$ must be retained.  In particular, the object
case has coefficient $(-1)^r$.  These are the intermediate signed
components; the final component also carries the output factor in
Equation~\eqref{eq:rh-higher-standard}.

For a unary component of form degree $p$ and internal degree
$q=|\phi|-p$, the normalization with trivial diagonal connection is
\begin{equation}\label{eq:unary-sign-calibration}
 \widetilde\RH_1(\phi^{(p,q)})(\sigma_p)
 =
 (-1)^q\int_{\Delta^p}\sigma_p^*\phi^{(p,q)},
 \qquad
 \RH_1(\phi^{(p,q)})(\sigma_p)
 =
 \int_{\Delta^p}\sigma_p^*\phi^{(p,q)}.
\end{equation}
The first formula follows from \cite[Definition 3.23]{AS12} after the
$J$-conjugation; the second follows by applying $\mathfrak G$ to the
internal-degree-$q$ output.  The binary calibration is then the
$n=2$ $A_\infty$ identity.

Let $(\psi^-)^V=\id_{\End(V)}\otimes\psi^-$ be the
endomorphism-valued $A_\infty$ de Rham map in the $-d$ convention.
The standard twisting formula for
$\widetilde\psi=\psi^-\circ J$ gives
\begin{equation}\label{eq:rh-higher}
 \begin{split}
 \widetilde\RH_n(\phi_n,\ldots,\phi_1)(\sigma_k)
 =\biggl[
 s^{-1}\sum_{a_0,\ldots,a_n\geq0}
 (\psi^-)^V_{n+\sum_i a_i}\bigl(
 &(s\alpha_n)^{\otimes a_n}\mid s\eta_n^J\mid
 \cdots\\[-2mm]
 &\cdots\mid s\eta_1^J\mid
 (s\alpha_0)^{\otimes a_0}
 \bigr)(\sigma_k)
 \biggr]_{n,0}.
 \end{split}
\end{equation}
If
\[
 q_{\mathrm{out}}
 =
 1-n+\sum_{i=1}^n|\phi_i|-k
\]
is the internal degree of this arity-$k$ output, then the final
Block--Smith component is
\begin{equation}\label{eq:rh-higher-standard}
 \RH_n(\phi_n,\ldots,\phi_1)(\sigma_k)
 =
 (-1)^{q_{\mathrm{out}}}
 \widetilde\RH_n(\phi_n,\ldots,\phi_1)(\sigma_k).
\end{equation}
Through Definition~\ref{def:chen} and the Igusa map, this is the
$[n,0]$ corner of the corresponding generalized-cone holonomy.  The
twisted $A_\infty$ identity has four types of terms:
\[
 D\phi_i,\qquad
 \phi_{i+1}\circ\phi_i,\qquad
 \RH_j\cup\RH_{n-j},\qquad
 D_{\Loc_{\mathrm{shHom}}}\RH_n.
\]
Their signed sum is zero.  This is exactly the $n$-th
$A_\infty$-functor identity: it holds first for
$\widetilde\RH$ by the twisting theorem and then for $\RH$ because
$\mathfrak G$ is a strict dg equivalence.  Thus
Equations~\eqref{eq:rh-higher} and
\eqref{eq:rh-higher-standard} supply a geometric explanation of the
higher components in Theorem~\ref{thm:integration}.
The unary cone calculation is \cite[Section 3.5]{BS14}; the complete
matrix and suspension calculation is
\cite[Definition 4.16 and Theorem 4.18]{AS12}.

\subsection{Quasi-full faithfulness}

The vertex complex of $\RH(E,\mathbb E)$ is
$(E_x,\mathbb E_x^0)$, so $\RH$ carries cohomologically regular
cohesive modules to cohomologically regular local systems.

\begin{lem}[Leading term of $\RH_1$]\label{lem:rh-leading}
Filter a source morphism by leafwise form degree and a target morphism
by simplicial arity.  Then $\RH_1$ preserves these decreasing
filtrations.  After taking cohomology with respect to the internal
differentials, put
\[
 H_i=H^\bullet(E_i,\mathbb E_i^0),
 \qquad
 L=\underline{\operatorname{Hom}}(H_1,H_2).
\]
The map on the resulting $E_1$-complexes is leafwise integration with
coefficients in $L$ and parallel transport for its induced flat
connection, after the two-ended target complex is identified as in
Lemma~\ref{lem:target-e1}.
\end{lem}

\begin{proof}
Write $\phi^{(p,q)}$ for the form-degree-$p$, internal-degree-$q$
component of a homogeneous source morphism, so $|\phi|=p+q$.  The
corner block of the generalized-cone
holonomy is a sum of Chen iterated integrals containing exactly one
copy of $\phi^{(p,q)}$; all other factors are positive-form
components of $\mathbb E_1$ and $\mathbb E_2$.  By
Definition~\ref{def:chen}, every nonempty word containing a zero-form
input vanishes.  Thus the form-degree-zero superconnection components
do not occur as nonzero Chen inputs; they remain present in the
internal differential and in the Stokes and bar identities.  Likewise,
all nonempty corner words containing $\phi^{(0,q)}$ vanish.  The
final arity-zero value is defined separately and is
$\phi^{(0,q)}$ on vertices by
Equation~\eqref{eq:categorical-unary-sign}.

Suppose now that $p\geq1$ and that a nonzero corner word contains $m$
diagonal inputs of form degrees $r_1,\ldots,r_m\geq1$.  Its
ordered-simplex fiber has dimension $m+1$, including the time assigned
to $\phi^{(p,q)}$.  The resulting path-space form therefore has degree
\[
 p+\sum_{i=1}^m r_i-(m+1).
\]
An arity-$k$ component is obtained by integration over the Igusa cube
$I^{k-1}$, so this degree must equal $k-1$.  Hence
\begin{equation}\label{eq:rh-leading-degree}
 k-p=\sum_{i=1}^m(r_i-1)\geq0.
\end{equation}
The internal degrees give the same equality: the insertions have
degrees $1-r_i$, the input has degree $|\phi|-p$, and the output has
degree $|\phi|-k$.  Thus form degree $p$ contributes only to arities
$k\geq p$, and an insertion of form degree at least two raises arity
strictly.

For $p\geq1$, equality $k=p$ occurs exactly when every diagonal input
is a connection one-form.  Including the empty words on either side,
the intermediate sum has a factor $(-1)^{a+b+q}$, and target grading
conjugation contributes the output factor $(-1)^q$.  Thus the complete
final equality sum along a path $\gamma$ is
\[
 \sum_{a,b\geq0}
 (-1)^{a+b}
 \operatorname{It}_\gamma
 \bigl((\mathbb E_2^1)^{\otimes b}\mid
       \phi^{(p,q)}\mid
       (\mathbb E_1^1)^{\otimes a}\bigr)
 =
 \int_0^1
 P_2(1,t)\,
 \bigl(\iota_{\dot\gamma(t)}\phi^{(p,q)}\bigr)
 P_1(t,0)\,dt,
\]
where $P_i$ denotes parallel transport for $\mathbb E_i^1$ and the
contraction is in the time direction of the path family; the remaining
$p-1$ arguments lie in the parameter cube.  This is parallel
transport in the Hom coefficient bundle.  Pulling it back along
$\Theta_{(p)}$, integrating over $I^{p-1}$, and applying
Lemma~\ref{lem:igusa-degree} gives
simplexwise covariant integration of $\phi^{(p,q)}$, in
agreement with Equation~\eqref{eq:unary-sign-calibration}.  For
$p=0$, only the separately defined arity-zero component
$\phi^{(0,q)}$ remains.

The equations of form degrees one and two in
$\mathbb E_i^2=0$ show that $\mathbb E_i^1$ is a chain connection and
that its curvature is internally null-homotopic.  It therefore induces
a flat connection on $H^\bullet(E_i,\mathbb E_i^0)$ and hence on $L$.
This identifies the leading map with the claimed coefficient
integration.  The same degree sorting is the unary generalized-cone
calculation in \cite[Sections 3.4--4.1]{BS14}; its Lie-algebroid
extension is \cite[Theorem 4.18]{AS12}.
\end{proof}

\begin{lem}[The two-ended target $E_1$-complex]
\label{lem:target-e1}
Let $(E_i,\mathbb E_i)$ be regular cohesive modules, let
\[
 H_i=H^\bullet(E_i,\mathbb E_i^0),
 \qquad
 L=\underline{\Hom}(H_1,H_2),
\]
and give $L$ the flat Hom connection.  After taking internal
cohomology, the arity-$p$ associated-graded target term consists of
plot-smooth maps
\[
 a_p(\sigma):
 H_1|_{\sigma(v_p)}
 \longrightarrow
 H_2|_{\sigma(v_0)}.
\]
It is naturally isomorphic, as a sheaf complex, to the one-sided
flat-coefficient cochain complex
$C^\bullet_{\mathrm{sh}}(\CF;L)$.  On a local raw representative the
isomorphism is
\begin{equation}\label{eq:two-ended-xi}
 \Xi_p(a)(\sigma)
 =
 a_p(\sigma)\circ
 \tau^1_{\sigma(v_0)\to\sigma(v_p)}
 \in L_{\sigma(v_0)}.
\end{equation}
\end{lem}

\begin{proof}
The arity-one holonomy of $\RH(E_i,\mathbb E_i)$ induces on $H_i$
the parallel transport $\tau^i$ of the flat cohomology connection.
The two-simplex Maurer--Cartan equation, after internal cohomology,
says that these arity-one maps compose strictly.  In particular their
restrictions to the edges of a simplex are invertible and independent
of the edge decomposition used for transport.

Equation~\eqref{eq:two-ended-xi} has the inverse
\[
 c(\sigma)\longmapsto
 c(\sigma)\circ
 \tau^1_{\sigma(v_p)\to\sigma(v_0)}.
\]
Both formulas are plot-smooth and commute with restriction, so they
give mutually inverse maps before and after sheafification.  It
remains to check the differential.  For a $(p+1)$-simplex, the
initial endpoint term in the bar differential postcomposes $a_p$ with
$\tau^2_{\sigma(v_1)\to\sigma(v_0)}$; under $\Xi$ this is transport
along the first edge in the Hom bundle,
\[
 T^L_{\sigma(v_1)\to\sigma(v_0)}(u)
 =
 \tau^2_{\sigma(v_1)\to\sigma(v_0)}
 \circ u\circ
 \tau^1_{\sigma(v_0)\to\sigma(v_1)}.
\]
The interior bar terms become the corresponding interior faces.  In
the terminal endpoint term, precomposition by
$\tau^1_{\sigma(v_{p+1})\to\sigma(v_p)}$ combines with
$\tau^1_{\sigma(v_0)\to\sigma(v_{p+1})}$ in
Equation~\eqref{eq:two-ended-xi}; strict composition turns it into the
last face.  With the suspended signs of the bar differential this
gives exactly
\[
 \Xi_{p+1}(d_1a)(\sigma)
 =
 T^L_{\sigma(v_1)\to\sigma(v_0)}
   \Xi_p(a)(\sigma\circ\partial_0)
 +\sum_{j=1}^{p+1}(-1)^j
   \Xi_p(a)(\sigma\circ\partial_j)
 =
 \delta_\nabla\Xi_p(a)(\sigma).
\]
Thus $\Xi$ is an isomorphism of the $E_1$-complexes.
\end{proof}

\begin{prop}[Quasi-full faithfulness]\label{prop:qff}
The restriction
\[
 \RH:\mathcal P_{\A}^{\mathrm{reg}}
 \longrightarrow
 \Loc_{\mathrm{shHom},\Ch_{\R}}^{\mathrm{sm,perf,reg}}
 \bigl(\Pi_{\infty}(\CF)\bigr)
\]
is $A_{\infty}$ quasi-fully faithful.
\end{prop}

\begin{proof}
For two source objects, apply $\RH_1$ on every open subset and
sheafify to obtain a morphism from their leafwise-form Hom sheaf to
the sheafified raw target Hom complex.  Filter the source by form
degree and the target by arity.  The source filtration is finite
because $\rk(\CF)<\infty$.  In a fixed target total degree only
finitely many arities occur by the amplitude estimate in
Proposition~\ref{prop:sheafified-target}.  Both filtrations are
therefore exhaustive, separated, and strongly convergent on every
stalk.

The $E_0$ differential is the commutator with the internal
differentials $\mathbb E_i^0$.  Regularity and bundle metrics give
smooth splittings, so internal cohomology is the vector bundle $L$ in
Lemma~\ref{lem:rh-leading}.  Lemma~\ref{lem:target-e1} identifies the
two-ended associated-graded target with the one-sided coefficient
complex, and Lemma~\ref{lem:rh-leading} identifies the induced map of
$E_1$-complexes with flat-coefficient integration.  The local part of
Theorem~\ref{thm:twisted-dr} makes the map a quasi-isomorphism on
every stalk; equivalently, the induced map on $E_2$ is an isomorphism.
Filtered comparison now shows that the full sheafified $\RH_1$ is a
stalkwise quasi-isomorphism.

Every source term is a sheaf of smooth bundle-valued leafwise forms
and hence is fine.  Every target term is fine by multiplication at the
first vertex, as in Theorem~\ref{thm:leafwise-dr}; this operation need
not commute with the differential.  Both total complexes are bounded
below.  The cone of the sheafified $\RH_1$ is thus a bounded-below,
sheafwise acyclic complex of fine sheaves, so its global-section
complex is acyclic.  Hence $\RH_1$ is a quasi-isomorphism on every
global morphism complex.  This is the sheaf-local analogue of
\cite[Proposition 4.3]{BS14}.
\end{proof}

\subsection{Prism transformations and local reconstruction}

Let $X_\bullet$ and $Y_\bullet$ be simplicial sets and let
\[
 \mathcal H:X_\bullet\times\Delta[1]\longrightarrow Y_\bullet
\]
be a simplicial homotopy from $f$ to $g$.  For
$0\leq j\leq n$, let $\lambda_j:[n+1]\to[1]$ take the value zero on
$0,\ldots,j$ and one on $j+1,\ldots,n+1$, and set
\[
 h_j(\sigma)=
 \mathcal H_{n+1}(s_j\sigma,\lambda_j),
 \qquad
 P_n(\sigma)=\sum_{j=0}^n(-1)^j h_j(\sigma).
\]
Thus $P$ is the oriented normalized prism from $f$ to $g$
\cite{May67}.

\begin{lem}[Strictly unital prism transformation]
\label{lem:prism-transformation}
Let $V$ be a strictly unital normalized Block--Smith
infinity-local system on $Y_\bullet$, with structure maps $F_m$.
Define
\begin{equation}\label{eq:unital-prism-components}
 q_0(x)=F_1(h_0x),
 \qquad
 q_n(\sigma)=
 \sum_{j=0}^n(-1)^jF_{n+1}(h_j\sigma)
 \quad(n\geq1).
\end{equation}
Then $q=\{q_n\}$ is a strictly unital normalized closed degree-zero
morphism
\[
 q:g^*V\longrightarrow f^*V.
\]
If $\mathcal H$ is plot-smooth and $V$ is smooth, then $q$ is a raw
plot-smooth morphism.
\end{lem}

\begin{proof}
In the unital normalized chain cocategory the prism satisfies
\begin{equation}\label{eq:prism-identities}
 \partial P+P\partial=g_*-f_*,
 \qquad
 \Delta_{\mathrm{AW}}P
 =
 (f_*\otimes P+P\otimes g_*)\Delta_{\mathrm{AW}}.
\end{equation}
The second formula is an identity of graded maps: since $P$ has chain
degree one,
\[
 (f_*\otimes P)(a\otimes b)
 =
 (-1)^{|a|}f_*(a)\otimes P(b).
\]
For example, if $\sigma=[v_0,v_1]$, $a_i=f(v_i)$, and $b_i=g(v_i)$,
then
\[
 P_1(\sigma)
 =
 [a_0,b_0,b_1]-[a_0,a_1,b_1].
\]
The middle Alexander--Whitney cuts have opposite signs, exactly as in
the graded formula.  In general, an Alexander--Whitney cut of a prism
simplex
\[
 h_j(\sigma)
 =
 [a_0,\ldots,a_j,b_j,\ldots,b_n]
\]
falls into exactly one of two cases.  If the cut is at position
$r\leq j$, its front factor is
$f_*(\sigma|[v_0,\ldots,v_r])$ and its back factor is
$h_{j-r}(\sigma|[v_r,\ldots,v_n])$.  The Koszul sign $(-1)^r$
from $f_*\otimes P$, multiplied by the prism sign
$(-1)^{j-r}$, is $(-1)^j$.  If $r\geq j+1$, put $s=r-1$.
The front factor is
$h_j(\sigma|[v_0,\ldots,v_s])$ and the back factor is
$g_*(\sigma|[v_s,\ldots,v_n])$, with prism sign $(-1)^j$.
These disjoint cases exhaust the Alexander--Whitney cuts and give the
two summands on the right of the second identity.  The first identity
is the usual prism boundary calculation.

Strict unitality is encoded in the unital, rather than reduced,
normalized chain cocategory.  A degenerate edge is zero in the
ordinary reduced normalized chain complex, whereas
$F_1(s_0y)=\id_{V_y}$.  Thus one must not interpret
\eqref{eq:unital-prism-components} as applying a reduced twisting
cochain to that quotient.  Equivalently, let
$\widetilde\tau_V$ be the unital extension of the twisting cochain to
the unital normalized dg cocategory in which degenerate edges
represent identity arrows and higher degenerate simplices vanish.
Then the component formula says
\[
 q=\widetilde\tau_VP.
\]
Here $\widetilde\tau_{V,m}=F_m$ in the signed convention fixed above;
no additional suspension factor is applied.
In particular, the constant homotopy gives the identity
transformation, including its arity-zero component.

Apply the unital twisting-cochain Maurer--Cartan equation to
$P_n\sigma$.
The first identity in \eqref{eq:prism-identities} supplies the
simplicial differential and the two endpoint terms; the second turns
the convolution term into
\[
 f^*\widetilde\tau_V\cup q
 +q\cup g^*\widetilde\tau_V
\]
with the required Koszul signs.  Their sum is precisely the
Block--Smith equation $Dq=0$.  Since $P_n\sigma$ has dimension $n+1$,
$q_n$ has internal degree $-n$, as required.  Its first and last
vertices are $f\sigma(v_0)$ and $g\sigma(v_n)$, so its endpoint type
is that of a morphism $g^*V\to f^*V$.

If $n\geq1$ and $\sigma$ is degenerate, $P_n\sigma$ is a sum of
degenerate simplices of dimension at least two, on which the strictly
unital structure maps vanish.  This proves normalization; arity zero
has no such vanishing condition.
If the homotopy is plot-smooth, every $h_j$ is plot-smooth; applying
the plot-smooth structure map of $V$ and taking the finite alternating
sum proves the final assertion.
\end{proof}

\begin{lem}[Raw Whitehead criterion]\label{lem:raw-whitehead}
Let $q:V\to W$ be a closed degree-zero morphism in
$\Loc_{\mathrm{raw}}$.  Suppose that its arity-zero component
\[
 q_0:(V,d_V)\longrightarrow(W,d_W)
\]
has a global smooth chain-homotopy inverse, with global smooth chain
homotopies.  Then $q$ is invertible in
$H^0(\Loc_{\mathrm{raw}})$.
\end{lem}

\begin{proof}
For every raw object $Z$, consider postcomposition and
precomposition,
\[
 q_*:\Hom_{\mathrm{raw}}(Z,V)\longrightarrow
      \Hom_{\mathrm{raw}}(Z,W),
 \qquad
 q^*:\Hom_{\mathrm{raw}}(W,Z)\longrightarrow
      \Hom_{\mathrm{raw}}(V,Z).
\]
Filter all four complexes by simplicial arity.  On the associated
graded for the internal differential, these maps are pointwise
postcomposition and precomposition by $q_0$.  A chosen smooth
chain-homotopy inverse to $q_0$, together with its smooth
homotopies, acts pointwise on every arity and supplies
chain-homotopy inverses on the associated graded.  All resulting
cochains remain plot-smooth.

For bounded finite-rank vertex complexes, only finitely many arities
occur in each total degree, by the amplitude estimate in
Proposition~\ref{prop:sheafified-target}.  The arity spectral
sequences are therefore strongly convergent, and filtered comparison
shows that $q_*$ and $q^*$ are quasi-isomorphisms.  Surjectivity of
$H^0(q_*)$ with $Z=W$ gives a right inverse
$[h_R]:W\to V$, and surjectivity of $H^0(q^*)$ with $Z=V$ gives a
left inverse $[h_L]:W\to V$.  In $H^0(\Loc_{\mathrm{raw}})$,
\[
 [h_L]
 =[h_L][q][h_R]
 =[h_R],
\]
so this common class is a two-sided inverse to $[q]$.  This is the
plot-smooth refinement of the arity-filtration argument in
\cite[Propositions 2.5 and 2.7]{BS14}.
\end{proof}

Write $\RH_{\mathrm{raw}}$ for the Block--Smith--Arias
Abad--Sch\"atz formulas before postcomposition with the germ functor
$\mathbf{Sh}$, and write
$\RH_{\mathrm{shHom}}=\mathbf{Sh}\circ\RH_{\mathrm{raw}}$.

\begin{thm}[Local higher Riemann--Hilbert correspondence]
\label{thm:local-rh}
Let $B=P\times T$ be a foliated box in which $P$ and $T$ are balls,
$P$ is star-shaped about $p_0$, and the plaques are
$P\times\{t\}$.  Then the raw and fixed-object sheafified-Hom
integration functors
\[
 \RH_{\mathrm{raw}}:
 \mathcal P_{\Omega^\bullet_{\CF}(B)}^{\mathrm{reg}}
 \longrightarrow
 \Loc_{\mathrm{raw},\Ch_{\R}}^{\mathrm{sm,perf,reg}}
 \bigl(\Pi_\infty(\CF|_B)\bigr)
\]
and
\[
 \RH_{\mathrm{shHom}}:
 \mathcal P_{\Omega^\bullet_{\CF}(B)}^{\mathrm{reg}}
 \longrightarrow
 \Loc_{\mathrm{shHom},\Ch_{\R}}^{\mathrm{sm,perf,reg}}
 \bigl(\Pi_\infty(\CF|_B)\bigr)
\]
are $A_\infty$ quasi-equivalences.
\end{thm}

\begin{proof}
Let $r:B\to T$ be projection and $i(t)=(p_0,t)$.  The contraction
\[
 H_s(p,t)=((1-s)p+sp_0,t)
\]
induces a simplicial homotopy from $\id$ to $ir$: for a monotone map
$\alpha:[n]\to[1]$, with affine realization
$|\alpha|:\Delta^n\to I$, put
\[
 \mathcal H_n(\sigma,\alpha)(u)
 =
 H_{|\alpha|(u)}(\sigma(u)).
\]
Every foliated simplex in $B$ has constant transverse coordinate.
It follows that the displayed simplex is again foliated.  Faces and
degeneracies commute with the construction because affine realization
is functorial, and the adjoint formula proves plot-smoothness.

Let $V$ be a target object.  Lemma~\ref{lem:prism-transformation}
gives a closed raw morphism
\begin{equation}\label{eq:local-prism-comparison}
 q_V:r^*i^*V\longrightarrow V,
 \qquad
 (q_V)_0(x)=F_1(h_0x),
 \qquad
 (q_V)_n(\sigma)=
 \sum_{j=0}^n(-1)^jF_{n+1}(h_j\sigma)
 \quad(n\geq1).
\end{equation}
For $x=(p,t)$ its arity-zero component is
\[
 (q_V)_0(x)=F_1(\gamma_x):
 V_{(p_0,t)}\longrightarrow V_x,
 \qquad
 \gamma_x(s)=H_s(x).
\]
Let $\bar\gamma_x(s)=\gamma_x(1-s)$.  The two foliated
$2$-simplices
\begin{align*}
 \Sigma_x(u_0,u_1,u_2)
 &=
 \bigl((u_0+u_2)p+u_1p_0,t\bigr),\\
 \bar\Sigma_x(u_0,u_1,u_2)
 &=
 \bigl((u_0+u_2)p_0+u_1p,t\bigr)
\end{align*}
have boundary paths
\[
 (\gamma_x,\bar\gamma_x,\operatorname{const}_x)
 \quad\text{and}\quad
 (\bar\gamma_x,\gamma_x,\operatorname{const}_{(p_0,t)}),
\]
respectively.  Their Maurer--Cartan equations and strict unitality
show that $F_1(\bar\gamma_x)$ is a two-sided chain-homotopy inverse to
$F_1(\gamma_x)$.  These maps and homotopies vary smoothly with $x$.

Put $K=i^*V$ and let
\[
 E_V=(r^*K,d_K,\nabla^{\mathrm{triv}})
\]
be the pullback complex with its trivial flat leafwise connection and
no higher connection components.  Since $T$ is a ball, the bounded
finite-rank graded bundle $K$ is trivializable, so $E_V$ is cohesive;
regularity is inherited from $V$.  Strict unitality and normalization
give the literal equality
\[
 \RH(E_V)=r^*i^*V.
\]
The unary transport is the identity and all normalized higher
transports vanish.  The construction above gives a global smooth
chain-homotopy inverse to $(q_V)_0$.  Therefore
Lemma~\ref{lem:raw-whitehead} makes
\eqref{eq:local-prism-comparison} an equivalence in the raw target.
Applying the strict dg germ functor $\mathbf{Sh}$ gives an equivalence
in the fixed-object sheafified-Hom target as well.

It remains only to note quasi-full faithfulness of the raw local
functor.  The coefficient contraction in the proof of
Theorem~\ref{thm:twisted-dr} is already a deformation retraction of
the raw coefficient complexes on $B$.  Applying the finite
form-degree/arity filtration and
Lemma~\ref{lem:rh-leading} before sheafification therefore proves that
the raw $\RH_1$ is a quasi-isomorphism.  The sheafified-Hom case is
Proposition~\ref{prop:qff}.  Hence both functors are quasi-fully
faithful and essentially surjective.
\end{proof}

\section{Effective algebraic descent}
\label{sec:cohesive-descent}

\subsection{Twisted leafwise cohesive modules}

For an open set $U\subset M$, write
\[
 \A_{\CF}(U)
 =
 \bigl(\Omega^\bullet_{\CF|_U}(U),d_{\CF}\bigr),
 \qquad
 \mathcal P_{\CF}(U)=\mathcal P_{\A_{\CF}(U)}.
\]
Write $\mathcal C_{\CF}(U)$ for the corresponding dg-category of
quasi-cohesive modules.  Restriction makes both
$U\mapsto\mathcal P_{\CF}(U)$ and
$U\mapsto\mathcal C_{\CF}(U)$ dg-presheaves.
Indeed, if a finite-projective module is represented by an idempotent
matrix, restriction applies entrywise to that idempotent; no flatness
of $C^\infty(U)\to C^\infty(V)$ is needed.

Fix a finite ordered open cover
$\mathfrak U=\{U_1,\ldots,U_N\}$ and put
$U_{i_0\ldots i_p}=U_{i_0}\cap\cdots\cap U_{i_p}$.
We recall the twisted model in the right-module convention of
\cite{BHW17,Wei23}.

\begin{defn}\label{def:twisted-cohesive}
A \emph{twisted leafwise cohesive module} on $\mathfrak U$ consists of
bounded graded finite-projective $C^\infty(U_i)$-modules
$E_i^\bullet$ and maps
\[
 a^{p,1-p}_{i_0\ldots i_p}:
 E_{i_p}^\bullet|_{U_{i_0\ldots i_p}}
 \otimes\A_{\CF}(U_{i_0\ldots i_p})
 \longrightarrow
 E_{i_0}^\bullet|_{U_{i_0\ldots i_p}}
 \otimes\A_{\CF}(U_{i_0\ldots i_p}),
 \qquad p\geq 0,
\]
of internal degree $1-p$.  The terms $a_i^{0,1}$ are flat
$\Z$-connections over $\A_{\CF}(U_i)$; the terms with $p\geq1$ are
$\A_{\CF}(U_{i_0\ldots i_p})$-linear.  They satisfy
\begin{equation}
 \delta a+a\cdot a=0,
 \label{eq:twisted-mc}
\end{equation}
and $a_{ii}^{1,0}$ is invertible up to homotopy.
\end{defn}

Here the product of homogeneous cochains is
\[
 (u\cdot v)_{i_0\ldots i_{p+r}}
 =
 (-1)^{qr}
 u^{p,q}_{i_0\ldots i_p}
 \circ
 v^{r,s}_{i_p\ldots i_{p+r}},
\]
and the endpoint-preserving \v{C}ech differential is
\[
 (\delta u)_{i_0\ldots i_{p+1}}
 =
 \sum_{k=1}^{p}(-1)^k
 u_{i_0\ldots\widehat{i_k}\ldots i_{p+1}}.
\]
Equation~\eqref{eq:twisted-mc} includes flatness of the local
$\Z$-connections and all higher cocycle identities.  If
$\mathcal E=(E_i,a)$ and $\mathcal F=(F_i,b)$, their morphism complex
is the total \v{C}ech complex of the local $\A_{\CF}$-linear Hom
complexes, with differential
\[
 D\phi=\delta\phi+b\cdot\phi-(-1)^{|\phi|}\phi\cdot a.
\]
We denote the resulting dg-category by
\[
 \operatorname{Tw}
 \bigl(M,\mathcal P_{\CF},\mathfrak U\bigr).
\]
The nondegeneracy condition ensures that this is homotopy-coherent
descent rather than an arbitrary Maurer--Cartan totalization.

Restriction of a global cohesive module defines the twisting functor
\[
 \mathcal T_{\mathfrak U}:
 \mathcal P_{\A}
 \longrightarrow
 \operatorname{Tw}
 \bigl(M,\mathcal P_{\CF},\mathfrak U\bigr):
\]
the local connections are the restrictions of $\mathbb E$, the
degree-$(1,0)$ transition maps are identities, and all higher
transition maps vanish.

\subsection{Sheafification and effectivity}

Let $\mathcal E=(E_i,a)$ be twisted.  For the open embedding
$j_I:U_I\hookrightarrow M$, use ordinary direct image,
\[
 (j_I)_*\mathcal G(V)=\mathcal G(V\cap U_I);
\]
this is not extension by zero and need not be locally free near the
boundary of $U_I$.  Define the sheafified total module by
\begin{equation}
 \mathcal S^n(\mathcal E)
 =
 \prod_{p+q=n}\prod_{i_0,\ldots,i_p}
 (j_{i_0\ldots i_p})_*
 \bigl(E_{i_0}^q|_{U_{i_0\ldots i_p}}\bigr)
 \label{eq:twisted-sheafification}
\end{equation}
and set
\[
 \mathbb S=\delta+a.
\]
We use Wei's product model.  The local internal degrees have a
common lower bound and the \v{C}ech degree is nonnegative, so the total
complex is bounded below; it need not be bounded in the \v{C}ech
direction.

If $\phi:\mathcal E\to\mathcal F$ is a homogeneous twisted morphism,
define $\mathcal S(\phi)$ on the product totalizations by the same
\v{C}ech convolution formula as the product above.  Thus
\[
 (\delta+b)\mathcal S(\phi)
 -(-1)^{|\phi|}\mathcal S(\phi)(\delta+a)
 =
 \mathcal S(D\phi).
\]
Consequently $\mathcal S$ is a dg functor on morphisms, not merely an
object construction.  Write
\[
 \mathbf S_{\mathfrak U}(\mathcal E)
 =
 \bigl(\Gamma(M,\mathcal S(\mathcal E)),\mathbb S\bigr)
\]
for its global quasi-cohesive totalization.

\begin{lem}\label{lem:sections-projective-tensor}
Let $\mathcal Q$ be a sheaf of $C^\infty_M$-modules and let
$V\to M$ be a finite-rank vector bundle with sheaf of sections
$\mathcal V$.  If $M$ is compact, the natural map
\[
 \Gamma(M,\mathcal Q)
 \otimes_{C^\infty(M)}\Gamma(M,V)
 \longrightarrow
 \Gamma\bigl(
 M,\mathcal Q\otimes_{C^\infty_M}\mathcal V
 \bigr)
\]
is an isomorphism.
\end{lem}

\begin{proof}
By Serre--Swan there is an idempotent
$e\in M_N(C^\infty(M))$ such that
$\Gamma(M,V)=eC^\infty(M)^N$ and
$\mathcal V=e(C^\infty_M)^N$.  Both sides of the displayed map
therefore identify naturally with
$e\,\Gamma(M,\mathcal Q)^N$.
\end{proof}

\begin{lem}\label{lem:twisted-quasi-cohesive}
If $M$ is compact, the pair
$\mathbf S_{\mathfrak U}(\mathcal E)$ is a quasi-cohesive module over
$\A$, and $\mathbf S_{\mathfrak U}$ is a dg functor to the dg-category
$\mathcal C_{\A}$ of quasi-cohesive $\A$-modules.
\end{lem}

\begin{proof}
Every positive \v{C}ech-degree component of $a$, as well as $\delta$,
is $\A$-linear.  The components $a_i^{0,1}$ obey the leafwise Leibniz
rule.  Consequently
\[
 \mathbb S(e\omega)
 =
 \mathbb S(e)\omega+(-1)^{|e|}e\,d_{\CF}\omega.
\]
A priori, the form-valued components of this operator take values in
\[
 \Gamma\bigl(
 M,\mathcal S(\mathcal E)
 \otimes_{C^\infty_M}\Omega^r_{\CF}
 \bigr).
\]
Lemma~\ref{lem:sections-projective-tensor}, applied to
$V=\extp^r\CF^\vee$, identifies this module with
\[
 \Gamma(M,\mathcal S(\mathcal E))
 \otimes_{C^\infty(M)}\Omega^r_{\CF}(M).
\]
Thus $\mathbb S$ is a $\Z$-connection over the global dga $\A$.
With the product sign in Definition~\ref{def:twisted-cohesive}, the
square of $\mathbb S$ is precisely
$\delta a+a\cdot a$ together with the local connection squares.
Equation~\eqref{eq:twisted-mc} therefore gives $\mathbb S^2=0$.
The total module need not be finite projective, which is why the result
is quasi-cohesive rather than cohesive.
\end{proof}

\begin{lem}[Twisting--sheafification adjunction]
\label{lem:twisting-sheafification-adjunction}
Assume that $M$ is compact.  In the quasi-cohesive enlargement,
restriction and product
totalization form an adjoint pair of dg functors
\[
 \mathcal T_{\mathfrak U}:
 \mathcal C_{\A}
 \rightleftarrows
 \operatorname{Tw}
 \bigl(M,\mathcal C_{\CF},\mathfrak U\bigr):
 \mathbf S_{\mathfrak U}.
\]
In particular, for $Q\in\mathcal C_{\A}$ and a twisted
quasi-cohesive object $\mathcal E$, there is a natural isomorphism of
complexes
\begin{equation}\label{eq:twisting-sheafification-adjunction}
 \Hom_{\operatorname{Tw}}
 \bigl(\mathcal T_{\mathfrak U}Q,\mathcal E\bigr)
 \cong
 \Hom_{\mathcal C_{\A}}
 \bigl(Q,\mathbf S_{\mathfrak U}\mathcal E\bigr).
\end{equation}
Denote the unit by
\[
 \epsilon_Q:Q\longrightarrow
 \mathbf S_{\mathfrak U}\mathcal T_{\mathfrak U}Q
\]
and the counit by
\[
 \eta_{\mathcal E}:
 \mathcal T_{\mathfrak U}\mathbf S_{\mathfrak U}\mathcal E
 \longrightarrow\mathcal E.
\]
\end{lem}

\begin{proof}
A morphism on the left of
Equation~\eqref{eq:twisting-sheafification-adjunction} is a compatible
family of \v{C}ech components.  Acting with those components on the
degree-zero restriction of a global section gives the morphism on the
right.  Conversely, projecting a map into the product totalization to
each \v{C}ech factor recovers the twisted components.  The two
operations are inverse, and the convolution identity preceding
Lemma~\ref{lem:sections-projective-tensor} shows that they commute
with the differentials and compositions.

This is the construction and adjunction of
\cite[Proposition 3.4]{Wei23}.  Its proof uses only restriction,
direct image, the \v{C}ech product, and the Leibniz rule for the local
dga; no Dolbeault-specific operation occurs.  It therefore applies
verbatim to the leafwise dga.  The identity morphisms under the
adjunction give the displayed unit and counit.
\end{proof}

\begin{lem}[Full unit comparison]\label{lem:full-descent-unit}
Assume that $M$ is compact.  For a global cohesive module $F$, the
adjunction unit
\[
 \epsilon_F:
 F\longrightarrow
 \mathbf S_{\mathfrak U}\mathcal T_{\mathfrak U}F
\]
is a quasi-isomorphism of quasi-cohesive modules.  More strongly, the
induced map of total $\A$-complexes
\[
 (F^\bullet\otimes_{C^\infty(M)}\A,\mathbb F)
 \longrightarrow
 \bigl(
  \mathbf S_{\mathfrak U}\mathcal T_{\mathfrak U}F
  \otimes_{C^\infty(M)}\A,\mathbb S
 \bigr)
\]
is a quasi-isomorphism.
\end{lem}

\begin{proof}
Filter both total complexes decreasingly by leafwise form degree.
The filtration is finite because
$\Omega^r_{\CF}=0$ for $r>\rk\CF$.  At form degree $r$, the associated
graded of the cone of $\epsilon_F$ is the augmented \v{C}ech complex
of the bounded complex of fine sheaves
\[
 \mathcal F^\bullet\otimes_{C^\infty_M}\Omega^r_{\CF}.
\]
A partition of unity subordinate to $\mathfrak U$ contracts this
augmented \v{C}ech complex in every internal degree.  Hence every
associated-graded cone is acyclic.  The finite filtered comparison
therefore makes the full cone acyclic.  Taking $r=0$ gives, in
particular, the quasi-isomorphism required in
$\mathcal C_{\A}$.  This also supplies the leafwise-form refinement
of the unit assertion in \cite[Proposition 3.4]{Wei23}.
\end{proof}

\begin{lem}\label{lem:leafwise-algebra-flat}
If $M$ is compact, then $\A$ is algebraically flat over
$C^\infty(M)$.
\end{lem}

\begin{proof}
For every $r$,
\[
 \Omega^r_{\CF}(M)
 =
 \Gamma\bigl(M,\extp^r\CF^\vee\bigr)
\]
is a finitely generated projective $C^\infty(M)$-module by
Serre--Swan.  Since it vanishes for $r>\rk\CF$, the underlying graded
module of $\A$ is a finite direct sum of projective modules and is
therefore flat.
\end{proof}

This lemma concerns the dga over its own degree-zero algebra.  It has
no bearing on the different question whether the abstract algebra
$C^\infty(P\times T)$ is flat over $C^\infty(T)$.

\begin{thm}[Effective descent for leafwise cohesive modules]
\label{thm:cohesive-descent}
Let $M$ be compact, let $\CF\subset TM$ be a regular foliation, and
let $\mathfrak U$ be a finite open cover.  The twisting functor is a
dg quasi-equivalence
\[
 \mathcal T_{\mathfrak U}:
 \mathcal P_{\A}
 \xrightarrow{\;\simeq_{\mathrm{qe}}\;}
 \operatorname{Tw}
 \bigl(M,\mathcal P_{\CF},\mathfrak U\bigr).
\]
\end{thm}

\begin{proof}
We adapt Wei's proof for the Dolbeault dga
\cite[Theorem 1.1]{Wei23}; we record why its effectivity step applies
to leafwise forms.

Let $\mathcal E=(E_i,a)$ be a twisted leafwise cohesive module and
first forget the positive leafwise-form components of $\mathbb S$.
The complex
\begin{equation}
 \bigl(\mathcal S(\mathcal E),\mathbb S^0\bigr)
 \label{eq:degree-zero-sheafification}
\end{equation}
is the sheafification of the underlying twisted perfect
$C^\infty_M$-complex.  Stalkwise, the degree-zero counit contracts the
simplex of cover indices containing the given point; this uses only
the twisted Maurer--Cartan equation and nondegeneracy.  Hence
Equation~\eqref{eq:degree-zero-sheafification} is locally
quasi-isomorphic to one of the bounded finite-rank complexes
$E_i^\bullet$ and is perfect.

Because $M$ is compact and the cover is finite, smooth
strict-perfect replacement gives a bounded complex of finite-rank
vector bundles $\mathcal E_0^\bullet$ and a quasi-isomorphism
\begin{equation}
 \phi_0:\mathcal E_0^\bullet
 \longrightarrow
 \bigl(\mathcal S(\mathcal E),\mathbb S^0\bigr).
 \label{eq:strict-perfect-replacement}
\end{equation}
It can be chosen so that the counit composite
\begin{equation}\label{eq:degree-zero-counit-equivalence}
 \eta^0_{\mathcal E}\circ
 \mathcal T_{\mathfrak U}(\phi_0):
 \mathcal T_{\mathfrak U}(\mathcal E_0^\bullet)
 \longrightarrow
 \operatorname{For}(\mathcal E)
\end{equation}
is a homotopy equivalence of twisted perfect
$C^\infty_M$-complexes.  This is
\cite[Proposition 4.3]{Wei23}: its proof combines smooth
strict-perfect replacement, the degree-zero counit, and the twisted
Whitehead criterion \cite[Proposition 3.3]{Wei23}.

All sheaves in Equation~\eqref{eq:strict-perfect-replacement} are
$C^\infty_M$-module sheaves and hence fine.  The cone is a bounded-below
acyclic complex of $\Gamma$-acyclic sheaves, so the acyclic-assembly
lemma, equivalently the hypercohomology spectral sequence, shows that
global sections preserve this quasi-isomorphism.  We obtain a bounded
finite-projective $C^\infty(M)$-complex
\[
 E_0^\bullet=\Gamma(M,\mathcal E_0^\bullet)
\]
and a $C^\infty(M)$-linear quasi-isomorphism
\begin{equation}
 \Gamma(M,\phi_0):E_0^\bullet
 \longrightarrow
 \bigl(\Gamma(M,\mathcal S(\mathcal E)),\mathbb S^0\bigr).
 \label{eq:global-degree-zero-model}
\end{equation}

Lemma~\ref{lem:leafwise-algebra-flat} now verifies the second
hypothesis of Block's quasi-representability theorem
\cite[Theorem 3.13]{Blo05}.  Applied to
$\mathbf S_{\mathfrak U}(\mathcal E)$ and the model
\eqref{eq:global-degree-zero-model}, it yields a cohesive module
$(E^\bullet,\mathbb E)\in\mathcal P_{\A}$ and a quasi-isomorphism
\[
 \phi:(E^\bullet,\mathbb E)
 \longrightarrow
 \mathbf S_{\mathfrak U}(\mathcal E)
\]
whose degree-zero forgetful image is
$\Gamma(M,\phi_0)$ under Serre--Swan.  This is the leafwise instance
of the construction in
\cite[Proposition 5.2]{Wei23}.

Use the counit from
Lemma~\ref{lem:twisting-sheafification-adjunction} to form the closed
morphism
\begin{equation}\label{eq:effective-descent-morphism}
 \Phi_{\mathcal E}
 =
 \eta_{\mathcal E}\circ\mathcal T_{\mathfrak U}(\phi):
 \mathcal T_{\mathfrak U}(E^\bullet,\mathbb E)
 \longrightarrow
 \mathcal E.
\end{equation}
After forgetting positive leafwise-form degree, this is
Equation~\eqref{eq:degree-zero-counit-equivalence}; in particular, its
local $(0,0)$ components are quasi-isomorphisms.  Since both source
and target are twisted cohesive objects, the twisted Whitehead
criterion \cite[Proposition 3.3]{Wei23} makes
$\Phi_{\mathcal E}$ a homotopy equivalence.  This proves essential
surjectivity and identifies the equivalence explicitly.

For global cohesive modules $E$ and $F$, the adjunction unit
\[
 \epsilon_F:
 F\longrightarrow
 \mathbf S_{\mathfrak U}\mathcal T_{\mathfrak U}(F)
\]
is the quasi-isomorphism of
Lemma~\ref{lem:full-descent-unit}.  Under
Equation~\eqref{eq:twisting-sheafification-adjunction}, the natural
restriction morphism
\[
 \operatorname{Hom}_{\mathcal P_{\A}}(E,F)
 \longrightarrow
 \operatorname{Hom}_{\operatorname{Tw}}
 \bigl(\mathcal T_{\mathfrak U}E,
       \mathcal T_{\mathfrak U}F\bigr)
\]
corresponds to postcomposition with $\epsilon_F$,
\[
 \Hom_{\mathcal P_{\A}}(E,F)
 \longrightarrow
 \Hom_{\mathcal C_{\A}}
 \bigl(E,\mathbf S_{\mathfrak U}\mathcal T_{\mathfrak U}F\bigr).
\]
By \cite[Proposition 3.9]{Blo05}, a quasi-isomorphism of
quasi-cohesive modules induces a quasi-isomorphism on the dg Yoneda
modules of cohesive objects.  Hence the restriction morphism is a
quasi-isomorphism.  Thus $\mathcal T_{\mathfrak U}$ is
quasi-fully faithful and essentially surjective.
\end{proof}

\begin{cor}\label{cor:cohesive-holim}
For a finite open cover of compact $M$, the global cohesive dg-category
is a model for the \v{C}ech homotopy limit of its restrictions:
\[
 \mathcal P_{\A}
 \simeq
 \operatorname*{holim}_{[p]\in\Delta}
 \prod_{i_0,\ldots,i_p}
 \mathcal P_{\A_{\CF}(U_{i_0\ldots i_p})}.
\]
\end{cor}

\begin{proof}
Block--Holstein--Wei identify the twisted dg-category with the
totalization, hence with the homotopy limit of the split \v{C}ech
diagram \cite[Theorem 4.1]{BHW17}.  Apply
Theorem~\ref{thm:cohesive-descent}.
\end{proof}

\begin{lem}[Dg nerves preserve the descent homotopy limit]
\label{lem:nerve-holim}
For a cosimplicial diagram $\mathcal C^\bullet$ of dg-categories,
there is a natural equivalence
\[
 N_{A_\infty}\left(
  \operatorname*{holim}_{[p]\in\Delta}\mathcal C^p
 \right)
 \simeq
 \lim_{[p]\in\Delta}
 N_{A_\infty}(\mathcal C^p)
\]
in $\operatorname{Cat}_\infty$.  Here the dg homotopy limit may be
computed after a Reedy-fibrant replacement; in particular, the
Block--Holstein--Wei totalization model has this property.
\end{lem}

\begin{proof}
On dg-categories, Faonte's $A_\infty$ nerve is naturally equivalent
to the dg nerve \cite[Section 2]{Fao17}.  The dg nerve is a right
Quillen functor from dg-categories with the Dwyer--Kan model structure
to simplicial sets with the categorical model structure
\cite[Proposition 1.3.1.20]{LurHA}.  Its right derived functor
therefore preserves homotopy limits.  A dg homotopy limit of a
fibrant replacement of $\mathcal C^\bullet$ is sent to the homotopy
limit of the corresponding dg nerves, and a homotopy limit in the
categorical model structure presents the limit in
$\operatorname{Cat}_\infty$.  Replacing the dg nerve by the naturally
equivalent $A_\infty$ nerve gives the displayed equivalence.
\end{proof}

\begin{lem}[Limits of fully faithful diagrams]
\label{lem:limit-fully-faithful}
Let
\[
 F^\bullet:\mathcal C^\bullet\longrightarrow\mathcal D^\bullet
\]
be a natural transformation of cosimplicial infinity-categories which
is fully faithful in every cosimplicial degree.  Then the induced
functor on limits is fully faithful.  Its essential image consists of
the limit objects whose component in every degree belongs to the
replete essential image of $F^p$.
\end{lem}

\begin{proof}
Mapping spaces in a limit of infinity-categories are limits of the
levelwise mapping spaces.  The first assertion therefore follows by
taking limits of equivalences.  For the second, replace
$\mathcal C^p$ by its replete essential image in $\mathcal D^p$.
Restriction functors preserve these full subcategories by naturality,
and the resulting levelwise equivalences induce an equivalence on
limits.
\end{proof}

\begin{lem}[Fiberwise detection of regularity]
\label{lem:fiberwise-regularity}
Evaluation at $x\in M$ defines a dg functor from leafwise cohesive
modules to finite complexes.  Consequently, an equivalence in the
$H^0$-category of cohesive modules evaluates to a chain-homotopy
equivalence.  A bounded finite-rank cohesive module which is locally
equivalent to regular cohesive modules is regular.
\end{lem}

\begin{proof}
Evaluation retains the zero-form component of a cohesive morphism.
The form-degree-zero part of its differential and composition is the
ordinary differential and composition of maps of fiber complexes.
It therefore sends a closed inverse and its two homotopies to a
chain-homotopy inverse on every fiber.

For the last assertion, local equivalence makes each
$\dim H^k(E_x)$ locally constant.  If
$e_k=\operatorname{rank} E^k$,
$r_k=\operatorname{rank}(\mathbb E^0:E^k\to E^{k+1})$, and
$h_k=\dim H^k(E_x)$, then
\[
 h_k=e_k-r_k-r_{k-1}.
\]
Starting in the lowest nonzero degree shows inductively that every
$r_k$ is locally constant.
\end{proof}

\begin{prop}[\v{C}ech full faithfulness]
\label{prop:shhom-cech-full-faithfulness}
Let $\mathfrak U=\{U_i\}_{i\in I}$ be a finite open cover of $M$,
and let $V,W$ be objects of
\[
 \Loc_{\mathrm{shHom},\Ch_{\R}}^{\mathrm{sm,perf}}
 \bigl(\Pi_\infty(\CF)\bigr)
\]
already defined globally on $M$.  Write
\[
 \check C_{\mathrm{full}}^{p,q}
 \bigl(\mathfrak U,\underline{\mathcal H}_{V,W}\bigr)
 =
 \prod_{(i_0,\ldots,i_p)\in I^{p+1}}
 \Gamma\bigl(
 U_{i_0\ldots i_p},
 \underline{\mathcal H}_{V,W}^{q}
 \bigr)
\]
for the full cosimplicial \v{C}ech--Hom bicomplex.  Then
$\underline{\mathcal H}_{V,W}^{\bullet}$ is bounded below and
termwise fine, and the augmentation
\[
 \Gamma\bigl(M,\underline{\mathcal H}_{V,W}^{\bullet}\bigr)
 \longrightarrow
 \operatorname{Tot}_{\Pi}
 \check C_{\mathrm{full}}^{\bullet,\bullet}
 \bigl(\mathfrak U,\underline{\mathcal H}_{V,W}\bigr)
\]
is a quasi-isomorphism.  The same holds for normalized cosimplicial
totalization.

Consequently, if $U_\bullet$ is the full \v{C}ech nerve of
$\mathfrak U$, restriction induces a fully faithful functor
\[
 N_{A_\infty}\left(
 \Loc_{\mathrm{shHom},\Ch_{\R}}^{\mathrm{sm,perf}}
 \bigl(\Pi_\infty(\CF)\bigr)\right)
 \longrightarrow
 \lim_{[p]\in\Delta}
 N_{A_\infty}\left(
 \Loc_{\mathrm{shHom},\Ch_{\R}}^{\mathrm{sm,perf}}
 \bigl(\Pi_\infty(\CF|_{U_p})\bigr)\right).
\]
Here the limit is taken in $\operatorname{Cat}_\infty$.
This is a mapping-descent statement for global objects; it makes no
claim of essential surjectivity or target object descent.
\end{prop}

\begin{proof}
The amplitude estimate in
Proposition~\ref{prop:sheafified-target} gives a lower bound for the
Hom complex.  In every fixed total degree, multiplication by a smooth
function at the first vertex gives the fineness endomorphisms used in
Theorem~\ref{thm:leafwise-dr}.

Choose a partition of unity $\{\rho_i\}_{i\in I}$ subordinate to the
cover.  In every fixed internal degree $q$, the usual formula
\[
 (hc)_{i_0\ldots i_{p-1}}
 =
 \sum_{j\in I}
 \Theta_{\rho_j}
 \bigl(c_{j i_0\ldots i_{p-1}}\bigr),
\]
with each summand extended by zero from
$U_{j i_0\ldots i_{p-1}}$, contracts the augmented cosimplicial
\v{C}ech row.  The extension is defined because the support of
$\Theta_{\rho_j}$ lies in $\supp(\rho_j)\subset U_j$.  The fineness
operators need not commute with the internal differential; only
rowwise \v{C}ech exactness is used.

If the Hom complex vanishes below $q_0$, then in total degree $n$
\[
 \operatorname{Tot}_{\Pi}^{n}
 =
 \prod_{p\geq0}
 \check C_{\mathrm{full}}^{p,n-p}
\]
has no terms for $p>n-q_0$.  Thus the apparently infinite product is
degreewise finite.  The spectral sequence obtained by first taking
\v{C}ech cohomology therefore converges and is concentrated in
\v{C}ech degree zero, where it is the global Hom complex.  This proves
the first assertion.  Standard cosimplicial normalization is a
deformation retract, proving the normalized statement.

Mapping spaces in a limit of infinity-categories are limits of the
levelwise mapping spaces.  Under the $A_\infty$ nerve, a
quasi-isomorphism of Hom complexes induces an equivalence of mapping
spaces, and the derived cosimplicial mapping complex is modeled by the
totalization above.  Restriction is therefore fully faithful.
\end{proof}

\begin{thm}[Compact higher Riemann--Hilbert correspondence]
\label{thm:compact-rh}
Let $M$ be compact and let $\CF\subset TM$ be a regular foliation.
Then
\[
 \RH_{\mathrm{shHom}}:
 \mathcal P_{\Omega^\bullet_{\CF}(M)}^{\mathrm{reg}}
 \longrightarrow
 \Loc_{\mathrm{shHom},\Ch_{\R}}^{\mathrm{sm,perf,reg}}
 \bigl(\Pi_\infty(\CF)\bigr)
\]
is an $A_\infty$ quasi-equivalence.
\end{thm}

\begin{proof}
Quasi-full faithfulness is Proposition~\ref{prop:qff}.  Choose a
finite cover $\mathfrak U=\{U_i\}$ by star-shaped foliated boxes.
For a global target object $V$, Theorem~\ref{thm:local-rh} gives
regular cohesive modules $E_i$ and raw equivalences
\[
 q_i:\RH(E_i)\xrightarrow{\;\sim\;}V|_{U_i}.
\]
Applying the strict germ functor $\mathbf{Sh}$ gives the corresponding
sheafified-Hom equivalences.

The regular source and target subcategories are full and replete.
For the source this follows from
Lemma~\ref{lem:fiberwise-regularity}; for the target, evaluation at a
vertex sends an equivalence and its inverse to a chain-homotopy
equivalence.  Restriction preserves regularity.  Their limits are
therefore the replete full subcategories of the ambient limits on
componentwise regular objects.

Apply the $A_\infty$ nerve to the source and target categories on the
full \v{C}ech diagram.  In every cosimplicial degree, the induced RH
functor is fully faithful by Proposition~\ref{prop:qff}, which applies
on arbitrary open intersections.  The restriction
$V|_{U_{i_0\ldots i_p}}$ lies in its replete essential image: it is
equivalent, by the restriction of $q_{i_0}$, to
\[
 \RH\bigl(E_{i_0}|_{U_{i_0\ldots i_p}}\bigr).
\]
Lemma~\ref{lem:limit-fully-faithful} therefore lifts the complete
homotopy-coherent \v{C}ech object of $V$ to the source \v{C}ech
limit, with componentwise regular local objects.  No assumption on the
shape of the intersections is used.

For the ambient source diagram,
Corollary~\ref{cor:cohesive-holim} and
Lemma~\ref{lem:nerve-holim} give the explicit equivalence
\[
 N_{A_\infty}(\mathcal P_{\A})
 \simeq
 \lim_{[p]\in\Delta}
 N_{A_\infty}\left(
  \prod_{i_0,\ldots,i_p}
  \mathcal P_{\A_{\CF}(U_{i_0\ldots i_p})}
 \right).
\]
It represents the lifted descent object by a global cohesive module
$E$.  Its restriction to each $U_i$ is equivalent to $E_i$, so
Lemma~\ref{lem:fiberwise-regularity} shows that $E$ is regular.

It remains to descend the equivalence to the two already global target
objects $\RH(E)$ and $V$.  Their Hom complex is
\[
 \Gamma\bigl(
 M,\underline{\mathcal H}^\bullet_{\RH(E),V}
 \bigr).
\]
Proposition~\ref{prop:shhom-cech-full-faithfulness} says precisely
that restriction to the full $\operatorname{Cat}_\infty$-valued
\v{C}ech limit is fully faithful on these global objects.  Hence the
equivalence between their restrictions lifts uniquely to an equivalence
\[
 \RH(E)\simeq V
\]
in the global $H^0$-category.  Thus $\RH_{\mathrm{shHom}}$ is
essentially surjective.
\end{proof}

\begin{cor}[Higher Riemann--Hilbert correspondence for proper submersions]
\label{cor:proper-submersion}
Let $\pi:M\to T$ be a proper submersion with compact base $T$, and put
$\CF=\ker(d\pi)$.  Then $M$ is compact and
Theorem~\ref{thm:compact-rh} gives an $A_\infty$ quasi-equivalence
\[
 \mathcal P_{\Omega^\bullet_{\CF}(M)}^{\mathrm{reg}}
 \simeq
 \Loc_{\mathrm{shHom},\Ch_{\R}}^{\mathrm{sm,perf,reg}}
 \bigl(\Pi_\infty(\CF)\bigr).
\]
\end{cor}

\begin{proof}
Properness and compactness of $T$ imply compactness of $M$, while a
submersion has a regular vertical foliation.  Apply
Theorem~\ref{thm:compact-rh}.  Ehresmann local triviality gives a
fiber-bundle interpretation of the local models, but the proof of the
correspondence needs no hypothesis beyond compactness and regularity.
\end{proof}

\begin{rem}[Scope of the descent theorem]
\label{rem:descent-scope}
Theorem~\ref{thm:cohesive-descent} uses finitely generated projective
modules, ordinary tensor products, and algebraic quasi-isomorphisms.
Maps between its finite-rank section modules are continuous for their
natural Fr\'echet topologies, but the proof does not establish closed
images, strict exactness, completed-tensor compatibility, or a
continuous strong deformation retract of
Equation~\eqref{eq:degree-zero-sheafification}.  It is therefore not a
descent theorem in Block's quasi-perfect Fr\'echet category.
\end{rem}

Together with the local prism reconstruction,
Theorem~\ref{thm:cohesive-descent} globalizes the inverse objects
needed in Theorem~\ref{thm:compact-rh}.  Its scope is algebraic source
descent; raw target-object descent, strict Fr\'echet descent, and the
Raw Hom Comparison of Conjecture~\ref{conj:raw-hom} remain separate
problems.

\section{Examples}\label{sec:examples}

The examples below exhibit smooth transverse parameters, a compact
foliation with no reasonable leaf space, and a nontrivial fiberwise
cohomology bundle.  All cohomology groups use real coefficients.

\begin{example}[The product foliation on $S^1\times T$]
\label{ex:product-circle}
Let $T$ be a compact smooth manifold and give
\[
 M=S^1\times T,\qquad S^1=\R/\Z,
\]
the vertical foliation $\CF=TS^1$.  In the coordinate $\theta$,
\[
 \Omega^0_{\CF}(M)=C^\infty(S^1\times T),\qquad
 \Omega^1_{\CF}(M)=C^\infty(S^1\times T)\,d\theta,
 \qquad
 d_{\CF}f=\frac{\partial f}{\partial\theta}\,d\theta.
\]
A leafwise closed function is independent of $\theta$, so
\[
 H^0_{\CF}(S^1\times T)\cong C^\infty(T).
\]
Fiber integration
\[
 I(g\,d\theta)(t)=\int_0^1g(\theta,t)\,d\theta
\]
vanishes on exact forms.  Conversely, if $I(g\,d\theta)=0$, then
\[
 f(\theta,t)=\int_0^\theta g(s,t)\,ds
\]
is smooth and periodic, and $d_{\CF}f=g\,d\theta$.  Hence
\[
 H^1_{\CF}(S^1\times T)\cong C^\infty(T),
 \qquad
 H^k_{\CF}(S^1\times T)=0\quad(k\geq2).
\]
More generally, if $L$ is compact, a fixed metric on $L$ and
fiberwise Hodge theory give
\[
 \Omega^\bullet_{\CF}(L\times T)
 \cong C^\infty\bigl(T,\Omega^\bullet(L)\bigr),
 \qquad
 H^k_{\CF}(L\times T)
 \cong C^\infty\bigl(T,H^k_{\mathrm{dR}}(L)\bigr)
\]
for the vertical foliation.

There is also an elementary description of the strict objects
concentrated in degree zero.  Let $E\to S^1\times T$ be a
finite-rank bundle with a flat partial connection along $TS^1$, and
put $K=E|_{\{0\}\times T}$.  Parallel transport on the interval gives
a smooth trivialization after pulling back to $[0,1]\times T$.  For
the positive loop $\ell_t(s)=(s,t)$, the last-to-first convention
gives a smooth bundle automorphism
\[
 A_t=F_1(\ell_t):K_t\longrightarrow K_t.
\]
Conversely, from a vector bundle $K\to T$ and a smooth bundle
automorphism $A:K\to K$, form
\[
 E_A=
 \operatorname{pr}_T^*K
 \big/
 \bigl((v,1)\sim(A^{-1}v,0)\bigr).
\]
The interval-direction trivial connection descends, and its
last-to-first transport around the positive loop is $A$.  Thus strict
flat leafwise bundles are equivalent to pairs $(K,A)$; a morphism
$u:(K,A)\to(K',A')$ is a bundle map satisfying $uA=A'u$.

The compact correspondence theorem extends this description from bundles to
globally bounded regular complexes with flat leafwise
superconnections and from strict monodromy to arbitrary plot-smooth
homotopy-coherent transport.  It identifies the full derived
morphism complexes, not merely the commuting maps between the
degree-zero monodromies.
\end{example}

\begin{example}[An irrational linear foliation]
\label{ex:irrational-torus}
Let $M=T^2=\R^2/\Z^2$, fix
$\alpha\in\R\setminus\Q$, and let $\CF$ be generated by
\[
 X=\frac{\partial}{\partial x}
   +\alpha\frac{\partial}{\partial y}.
\]
Every leaf is dense and diffeomorphic to $\R$.  The restriction
$\vartheta=dx|_{\CF}$ satisfies $\vartheta(X)=1$, so
\[
 H^0_{\CF}(T^2)=\ker X,\qquad
 H^1_{\CF}(T^2)
 =
 C^\infty(T^2)/X\bigl(C^\infty(T^2)\bigr).
\]
For the complex Fourier modes
\[
 e_{m,n}(x,y)=e^{2\pi i(mx+ny)}
\]
one has
\[
 X(e_{m,n})=2\pi i(m+\alpha n)e_{m,n}.
\]
Irrationality implies that the only invariant smooth functions are
constant, and therefore
\[
 H^0_{\CF}(T^2)\cong\R.
\]

Suppose first that $\alpha$ is Diophantine: there are
$C>0$ and $\nu\geq0$ such that
\[
 |m+\alpha n|
 \geq C(1+|m|+|n|)^{-\nu}
\]
for every nonzero $(m,n)\in\Z^2$.  If a smooth zero-average function
$f$ has Fourier coefficients $\widehat f_{m,n}$, define
\[
 \widehat g_{m,n}
 =
 \frac{\widehat f_{m,n}}{2\pi i(m+\alpha n)}
 \quad ((m,n)\neq(0,0)),
 \qquad
 \widehat g_{0,0}=0.
\]
The coefficients of $f$ are rapidly decreasing, and the displayed
lower bound shows that those of $g$ remain rapidly decreasing.
The reality symmetry is preserved, so $g$ is real when $f$ is real,
and $Xg=f$.  The average therefore induces an isomorphism
\[
 H^1_{\CF}(T^2)\cong\R
\]
in the Diophantine case.

Now suppose that $\alpha$ is Liouville.  Choose distinct
$k_j=(m_j,n_j)\in\Z^2\setminus\{0\}$ with
\[
 2^j\leq |k_j|:=|m_j|+|n_j|<|k_{j+1}|,\qquad
 \delta_j:=m_j+\alpha n_j\neq0,\qquad
 |\delta_j|\leq(1+|k_j|)^{-j}.
\]
Such a sequence is obtained by taking successively better Liouville
approximations with increasing denominators.  In particular, the mode
pairs $\{\pm k_j\}$ are pairwise disjoint.  Define
\[
 f=
 \sum_{j=1}^{\infty}
 2\pi i\delta_j
 \bigl(e_{m_j,n_j}-e_{-m_j,-n_j}\bigr).
\]
The paired coefficients are complex conjugates, so $f$ is real.  An
order-$r$ derivative of the $j$th pair is bounded by a constant times
\[
 |k_j|^r|\delta_j|
 \leq |k_j|^{r-j}
 \leq 2^{-j(j-r)}
 \qquad(j>r).
\]
Thus the Fourier series and all its derivatives converge uniformly,
and $f\in C^\infty(T^2)$.

If $Xg=f$, comparison of Fourier coefficients would give
\[
 \widehat g_{m_j,n_j}
 =
 \widehat g_{-m_j,-n_j}
 =1
\]
for every $j$, contradicting rapid decrease.  On the other hand, the
real trigonometric polynomials
\[
 g_N=\sum_{j=1}^N
 \bigl(e_{m_j,n_j}+e_{-m_j,-n_j}\bigr)
\]
satisfy $Xg_N=f_N$, where $f_N\to f$ in the Fr\'echet topology of
$C^\infty(T^2)$.  Hence $X(C^\infty(T^2))$ is not closed and the
natural quotient topology on $H^1_{\CF}(T^2)$ is non-Hausdorff.
This does not compute the full Liouville quotient.

Although this foliation has no smooth Hausdorff leaf space, $T^2$ is
compact, so Theorem~\ref{thm:compact-rh} still gives
\[
 \mathcal P_{\Omega^\bullet_{\CF}(T^2)}^{\mathrm{reg}}
 \simeq
 \Loc_{\mathrm{shHom},\Ch_{\R}}^{\mathrm{sm,perf,reg}}
 \bigl(\Pi_\infty(\CF)\bigr).
\]
The scalar calculation shows that the leafwise side can retain
arithmetic information invisible in a naive leaf-space description;
we make no elementary classification of all target objects.
\end{example}

\begin{example}[Mapping tori and cohomology bundles]
\label{ex:mapping-torus}
Let $L$ be compact, let $\varphi:L\to L$ be a diffeomorphism, and
form
\[
 M_\varphi
 =
 (L\times[0,1])/
 \bigl((x,1)\sim(\varphi(x),0)\bigr).
\]
The projection $\pi:M_\varphi\to S^1$ is a proper submersion; give
$M_\varphi$ the vertical foliation $\CF=\ker(d\pi)$.  A vertical form
is a smooth family $\omega_t\in\Omega^\bullet(L)$ satisfying
\[
 \omega_1=\varphi^*\omega_0
\]
and the corresponding compatibility for all endpoint derivatives.
Passing to cohomology gives a vector bundle
$\mathcal H_\varphi^k\to S^1$ with fiber
$H^k_{\mathrm{dR}}(L)$ and gluing convention
\[
 (a,1)\sim\bigl((\varphi^{-1})^*a,0\bigr).
\]
The explicit equation $a_1=\varphi^*a_0$ avoids any ambiguity about
the orientation convention for monodromy.

A chosen vertical metric identifies this cohomology bundle smoothly
with the corresponding bundle of fiberwise harmonic forms.
Parameter-dependent elliptic theory gives smooth harmonic
projections and Green operators, and hence
\[
 H^k_{\CF}(M_\varphi)
 \cong
 \Gamma^\infty(S^1,\mathcal H_\varphi^k).
\]
If $L=S^2$ and $\varphi$ reverses orientation, then
$\varphi^*=-1$ on $H^2_{\mathrm{dR}}(S^2)$.  The bundle
$\mathcal H_\varphi^2$ is therefore the M\"obius real line bundle, and
\[
 H^2_{\CF}(M_\varphi)
 \cong
 \Gamma^\infty(S^1,\mathcal H_\varphi^2).
\]
Corollary~\ref{cor:proper-submersion} identifies regular vertical
cohesive modules on $M_\varphi$ with regular transversely smooth
infinity-local systems.  This example shows that the correspondence
retains nontrivial transverse twisting even though the foliation is
locally a product.
\end{example}

\section{Topological enhancements and strict reconstruction}
\label{sec:topological-enhancements}

This section refines the algebraic theory under additional topological
hypotheses.  It gives sufficient conditions for compatible nuclear
Fr\'echet models and for the preservation of strict exact sequences
under completed tensor products.  The de Rham comparison,
quasi-full-faithfulness, local prism reconstruction, compact descent,
and hypersheafification results are independent of these hypotheses.
Throughout this section, ``exact'' refers to the explicitly stated
quasi-abelian exact structure.

\subsection{Fr\'echet modules and completed scalar extension}

Let $U=P\times T$ be a product foliation chart, with $P$ a nonempty
coordinate ball and plaques $P\times\{t\}$.  Put
\[
 A=C^\infty(T),\qquad
 E=C^\infty(P),\qquad
 B=C^\infty(P\times T).
\]
These spaces carry their standard nuclear Fr\'echet topologies.  Let
$\mathsf{FrMod}(A)$ be the category of unital Fr\'echet $A$-modules
with jointly continuous action and continuous module maps.
This is the relative Fr\'echet setting used for Block's quasi-perfect
category \cite[Section 2.2]{Blo06II}; here we make the completed
balanced tensor product and the strict exact structure explicit.

A short exact sequence in this category is \emph{strict} if the first
map is a topological embedding onto a closed subspace and the last map
identifies the target topologically with the quotient.  It is
\emph{admissibly split} if it also has a continuous linear splitting,
not necessarily an $A$-linear one.  A complex is \emph{strictly exact}
if every sequence
\[
 0\longrightarrow Z^n
 \longrightarrow C^n
 \xrightarrow{\,d\,}Z^{n+1}
 \longrightarrow0
\]
is strict, and a cochain map is a \emph{strict quasi-isomorphism} if
its mapping cone is strictly exact.  For a right Fr\'echet $A$-module
$M$ and a left Fr\'echet $A$-module $N$, set
\[
 M\mathbin{\widehat\otimes_{A,\pi}}N
 =
 \frac{M\widehat\otimes_\pi N}
 {\overline{\operatorname{span}}\{
 ma\otimes n-m\otimes an\}}.
\]
The quotient is again Fr\'echet.  This is the only completed relative
tensor product used below.

\begin{thm}[Local strict topological flatness]
\label{thm:top-flat}
With the notation above:
\begin{enumerate}
 \item multiplication induces an isomorphism of nuclear Fr\'echet
 algebras
 \[
  A\widehat\otimes_\pi E
  \xrightarrow{\;\cong\;}B,
  \qquad
  a\otimes e\longmapsto
  \bigl((p,t)\mapsto a(t)e(p)\bigr);
 \]
 \item for every $M\in\mathsf{FrMod}(A)$ there is a natural
 topological $B$-module isomorphism
 \[
  \beta_M:
  M\mathbin{\widehat\otimes_{A,\pi}}B
  \xrightarrow{\;\cong\;}
  M\widehat\otimes_\pi E,
  \qquad
  m\widehat\otimes_A(a\otimes e)\longmapsto ma\otimes e;
 \]
 \item completed extension of scalars
 \[
  -\mathbin{\widehat\otimes_{A,\pi}}B
 \]
 preserves strict short exact sequences of Fr\'echet modules,
 strictly exact complexes, and strict quasi-isomorphisms;
 \item if $p_0\in P$, evaluation at $p_0$ splits the unit of scalar
 extension.  In particular, this extension functor is faithful and
 conservative.
\end{enumerate}
\end{thm}

\begin{proof}
The smooth kernel theorem gives the first isomorphism
\cite[Theorem 51.6]{Tre67}; see also \cite{Gro55}.  It agrees with
multiplication on the dense subalgebra $A\otimes E$, hence is an
isomorphism of topological algebras.

Associativity of completed projective tensor products gives a
continuous map
\[
 M\widehat\otimes_\pi(A\widehat\otimes_\pi E)
 \longrightarrow M\widehat\otimes_\pi E,
 \qquad
 m\otimes a\otimes e\longmapsto ma\otimes e.
\]
It annihilates the closed balanced subspace and hence induces
$\beta_M$.  Conversely,
\[
 m\otimes e\longmapsto
 m\widehat\otimes_A(1\otimes e)
\]
induces a continuous map in the opposite direction.  The two maps are
inverse on elementary balanced tensors and therefore everywhere.

If
\[
 0\longrightarrow M'\longrightarrow M\longrightarrow M''
 \longrightarrow0
\]
is strict exact, scalar extension identifies it with the sequence
obtained by applying $-\widehat\otimes_\pi E$.  A nuclear Fr\'echet
space is exact on strict short exact sequences of Fr\'echet spaces for
the completed projective tensor product: projective tensor product
preserves strict quotients, injective tensor product preserves strict
embeddings, and the two completed tensor products agree when one factor
is nuclear
\cite[Propositions 43.7, 43.9 and Theorem 50.1(f)]{Tre67};
see also \cite[Chapter I, Section 1, Proposition 3 and Chapter II,
Section 3, Proposition 10]{Gro55}.  This proves strict exactness.  The
assertions for complexes follow by applying the same result to strict
cycles, boundaries, and mapping cones.  Admissibly split exactness,
when needed, follows already from preservation of finite topological
direct sums.

Finally,
\[
 \rho_{p_0}:M\widehat\otimes_\pi E\longrightarrow M,
 \qquad
 m\otimes e\longmapsto e(p_0)m,
\]
is a continuous $A$-linear retraction of the unit
$\eta_M:m\mapsto m\otimes1$.  Write
$F=-\widehat\otimes_{A,\pi}B$.
If scalar extension sends $f:M\to N$ to an isomorphism $F(f)$, then
\[
 \rho_M\circ F(f)^{-1}\circ\eta_N:N\longrightarrow M
\]
is inverse to $f$ by naturality of $\eta$ and $\rho$.  The same split
unit shows that $F(f)=0$ implies $f=0$.  Thus extension is conservative
and faithful.
\end{proof}

On this chart,
\[
 \Omega^r_{\CF}(U)
 \cong
 B\otimes\extp^r(\R^{\dim P})^\vee.
\]
For every Fr\'echet $A$-module $K$, completed balancing gives the
associativity isomorphism
\begin{align}
 \bigl(K\mathbin{\widehat\otimes_{A,\pi}}B\bigr)
 \mathbin{\widehat\otimes_{B,\pi}}\Omega^r_{\CF}(U)
 &\xrightarrow{\;\cong\;}
 K\mathbin{\widehat\otimes_{A,\pi}}\Omega^r_{\CF}(U),
 \notag\\
 (k\widehat\otimes_{A,\pi} b)
 \widehat\otimes_{B,\pi}\omega
 &\longmapsto k\widehat\otimes_{A,\pi} b\omega.
 \label{eq:relative-associativity}
\end{align}
Its inverse sends $k\widehat\otimes_{A,\pi}\omega$ to
$(k\widehat\otimes_{A,\pi}1)
\widehat\otimes_{B,\pi}\omega$.  Thus local finite models may be
extended first to smooth functions and then to leafwise forms without
changing the result.

\begin{rem}[Topological versus algebraic flatness]
\label{rem:flatness}
Theorem~\ref{thm:top-flat} is a theorem about completed projective
tensor products and the strict exact category of topological modules.
Its conclusion does not extend to arbitrary algebraic modules,
discontinuous maps, or algebraically exact sequences with nonclosed
images.  In particular, it does not imply that the underlying abstract
ring $C^\infty(P\times T)$ is flat over the underlying abstract ring
$C^\infty(T)$.  No assertion about ordinary algebraic flatness is made
here; the proof above controls only completed tensor products and
strict exact sequences.
\end{rem}

\subsection{Scope of local strict flatness}

Let $\mathbf R$ be the unit infinity-local system.  For a target object
$V$, let $\underline C_V^\bullet$ be the sheafification of the raw
unit-morphism presheaf
\[
 U\longmapsto
 \mathcal H^\bullet_{\mathbf R|_U,V|_U}(U),
\qquad
 C_V^\bullet(U)
 =
 \Gamma(U,\underline C_V^\bullet|_U).
\]
This is algebraically a complex of
$\underline{\R}_{\CF}$-modules.  The functional diffeology on the
spaces of simplices does not, by itself, put a canonical Fr\'echet
topology on this whole morphism complex.  Hence an expression such as
\[
 C_V^\bullet(U)
 \widehat\otimes_{C^\infty(T)}C^\infty(P\times T)
\]
is not defined until a suitable local topological model has been
constructed.

Nor does vertexwise perfectness imply that $C_V^\bullet(U)$ is perfect
over $C^\infty(T)$: the latter contains all simplicial degrees and the
full bar-convolution differential.  Theorem~\ref{thm:top-flat}
preserves a strict finite resolution once one exists, but it does not
produce that resolution or the coherent local cohesive descent datum,
and it does not prove that the evaluation counit is a quasi-isomorphism.
If such bounded finite-rank local cohesive objects and twisted
transition data have already been constructed,
Theorem~\ref{thm:cohesive-descent} supplies their algebraic
globalization.  It still does not provide the continuous strong
deformation retract used by the explicit strict model below.

\begin{defn}[Strict topological enhancement]
\label{def:strict-model}
A smooth perfect regular infinity-local system $V$ admits a
\emph{strict topological enhancement} if it is equipped with the
following data.
\begin{enumerate}
 \item \emph{Strict local finite models.}  On a finite product-chart
 cover $U_i=P_i\times T_i$, the small morphism complex
 $C_V^\bullet(U_i)$ has a restriction-compatible Fr\'echet
 $C^\infty(T_i)$-module model with a continuous strong deformation
 retract onto a bounded complex $K_i^\bullet$ of finitely generated
 projective $C^\infty(T_i)$-modules.
 \item \emph{Completed descent and strictification.}  The complexes
 \[
  K_i^\bullet
  \mathbin{\widehat\otimes_{C^\infty(T_i),\pi}}
  C^\infty(P_i\times T_i)
 \]
 and their leafwise-form extensions from
 Equation~\eqref{eq:relative-associativity} descend compatibly to a
 global quasi-cohesive sheaf complex $\mathcal Q_V$.  Its chartwise
 completed models are bounded projective Fr\'echet modules with
 continuous $\Z$-connections, hence are quasi-perfect in Block's
 topological sense.  Its form-degree zero part admits a continuous
 $C^\infty_M$-linear \emph{normalized} strong deformation retract onto
 a bounded finite-rank vector-bundle complex
 $(E^\bullet,\mathbb E^0)$.  Writing the differential on the
 form-degree zero part as $\mathbb F^0$, the retract consists of
 continuous maps $i,p,h$ with $|h|=-1$ and
 \[
  pi=\id,\qquad
  ip-\id=\mathbb F^0h+h\mathbb F^0,\qquad
  ph=hi=h^2=0.
 \]
 \item \emph{Integration and evaluation.}  In the chosen completed
 models, the higher-holonomy series for $\mathcal Q_V$ converge and
 glue to an infinity-local system, denoted
 $\operatorname{Int}(\mathcal Q_V)$.  The same convergence and gluing
 hold for the perturbed comparison maps obtained from the normalized
 retract in (2), so their integrations give mutually inverse
 equivalences between $\RH(E,\mathbb E)$ and
 $\operatorname{Int}(\mathcal Q_V)$ in
 $H^0(\Loc_{\mathrm{shHom}})$.  The chartwise radial
 evaluation maps are compatible on overlaps and induce a natural
 closed degree-zero morphism
 \[
  \varepsilon_V:
  \operatorname{Int}(\mathcal Q_V)\longrightarrow V.
 \]
 Its globally defined arity-zero component is required to be a
 quasi-isomorphism at every point.
\end{enumerate}
\end{defn}

\begin{prop}[Explicit strict reconstruction]
\label{prop:strict-reconstruction}
If $V$ admits a strict topological enhancement, the data of
Definition~\ref{def:strict-model} produce a regular cohesive module
$(E,\mathbb E)$ and a specified equivalence
\[
 \RH(E,\mathbb E)\xrightarrow{\;\sim\;}V
\]
in the fixed-object sheafified-Hom target.
\end{prop}

\begin{proof}
Theorem~\ref{thm:top-flat} preserves the strict
local finite models under completed scalar extension, and
Equation~\eqref{eq:relative-associativity} identifies their
leafwise-form extensions.  The second part of the enhancement supplies
a global strong deformation retract
\[
 \begin{tikzcd}
 (E^\bullet,\mathbb E^0)
 \arrow[r,shift left=1.1ex,"i"]
 &
 (\mathcal Q_V^{\bullet,0},\mathbb F^0)
 \arrow[l,shift left=1.1ex,"p"]
 \end{tikzcd}
 \qquad
 pi=\id,\quad
 ip-\id=\mathbb F^0h+h\mathbb F^0,\quad
 ph=hi=h^2=0.
\]

Write the flat differential of $\mathcal Q_V$ as
$\mathbb F^0+\delta$, where $\delta$ raises leafwise form degree.
Extend $i,p,h$ $C^\infty_M$-linearly to leafwise forms.  Because
$\Omega^r_{\CF}=0$ above $r=\rk\CF$, both $h\delta$ and $\delta h$ are
nilpotent.  Thus the perturbation is small in the sense of
\cite[Section 2]{Cra04}, and all of the following series are finite:
\[
 A_\delta
 =
 (1-\delta h)^{-1}\delta
 =
 \delta(1-h\delta)^{-1},
 \qquad
 \mathbb E
 =
 \mathbb E^0+pA_\delta i,
\]
with perturbed comparison data
\[
 i_\delta=i+hA_\delta i,\qquad
 p_\delta=p+pA_\delta h,\qquad
 h_\delta=h+hA_\delta h.
\]
The basic perturbation lemma gives
\[
 p_\delta i_\delta=\id,\qquad
 i_\delta p_\delta-\id
 =
 (\mathbb F^0+\delta)h_\delta
 +h_\delta(\mathbb F^0+\delta),
\]
as well as $\mathbb E^2=0$ and the normalized side conditions for the
perturbed retract.  Since $i,p,h$ are $C^\infty_M$-linear, the scalar
term in the Leibniz rule is unchanged; thus $\mathbb E$ is a flat
$\Z$-connection, not merely an $\R$-linear differential.  This is the
finite form-degree version of the transfer used in
\cite[Theorem 3.13]{Blo05} and in the proof of essential
surjectivity \cite[Theorem 4.5]{BS14}.

The third part of the enhancement applies higher holonomy to
$i_\delta,p_\delta,h_\delta$ and gives
\[
 \RH(E,\mathbb E)
 \xrightarrow[\simeq]{\operatorname{Int}(i_\delta)}
 \operatorname{Int}(\mathcal Q_V)
 \xrightarrow[\simeq]{\varepsilon_V}V.
\]
The first arrow is an equivalence by the integrated perturbed retract.
The second is an equivalence by
Proposition~\ref{prop:target-whitehead}, applied to the closed
morphism and its vertexwise quasi-isomorphic arity-zero component.
At vertices this comparison identifies
$H^\bullet(E_x,\mathbb E_x^0)$ with $H^\bullet(V_x)$, so $E$ is
cohomologically regular.  Thus it belongs to
$\mathcal P_{\A}^{\mathrm{reg}}$.
\end{proof}

\begin{rem}
The additional hypotheses of
Definition~\ref{def:strict-model} are independent of
Theorem~\ref{thm:compact-rh}, whose reconstruction uses local prisms
and algebraic dg descent.  Under these hypotheses,
Proposition~\ref{prop:strict-reconstruction} supplies an explicit
continuous finite model suitable for comparison with Block's
quasi-perfect Fr\'echet category or with the raw global target.  The
local strict-flatness theorem alone does not construct a strict
enhancement.
\end{rem}

Accordingly, the results of this section do not compare global raw
morphism complexes with their sheafified-Hom counterparts.

\begin{rem}
When $\CF=TM$, the compact correspondence specializes to the
Block--Smith quasi-equivalence \cite[Theorem 4.1]{BS14}.  When
$\CF=0$, both regular categories reduce directly to bounded complexes of
finite-rank smooth vector bundles with smooth transverse variation.
The smoothness condition in Definition~\ref{def:smooth-local-system}
excludes discontinuous families.
\end{rem}

\section{Localization and infinity-categorical formulation}

\subsection{Vertexwise quasi-isomorphisms and localization}

Let $W_{\mathrm{coh}}$ be the class of closed degree-zero morphisms
\[
 f:(E,\mathbb E)\longrightarrow(F,\mathbb F)
\]
in $\mathcal P_{\A}$ for which
\[
 f_x^0:(E_x,\mathbb E_x^0)\longrightarrow
 (F_x,\mathbb F_x^0)
\]
is a quasi-isomorphism for every $x\in M$.

\begin{prop}\label{prop:coh-localization}
Every morphism in $W_{\mathrm{coh}}$ is a homotopy equivalence in the
cohesive dg-category.
\end{prop}

\begin{proof}
The form-degree-zero part of $\operatorname{Cone}(f)$ is a pointwise
acyclic bounded complex of finite-rank vector bundles.  Exactness and
boundedness force the ranks of its differentials to be locally
constant.  After choosing bundle metrics, its fiberwise Hodge
Laplacian
\[
 \Delta=d_0d_0^*+d_0^*d_0
\]
is invertible and has a smooth inverse.  Hence
$h_0=d_0^*\Delta^{-1}$ contracts the degree-zero cone.

The remaining components of the cone $\Z$-connection raise leafwise
form degree.  Extending $h_0$ to leafwise forms and applying the finite
basic perturbation formula, which terminates above form degree
$\rk\CF$, contracts the total cone.  Thus $f$ is a homotopy
equivalence.  This is the foliated specialization of the cohesive
module criterion in \cite[Proposition 2.9]{Blo05}.
\end{proof}

Let $W_{\mathrm{shHom}}$ consist of the closed degree-zero morphisms
in $\Loc_{\mathrm{shHom}}$ whose arity-zero components are
quasi-isomorphisms at every point.  The raw arity-zero projection has
values in the sheaf of smooth bundle maps, which is already a sheaf;
its universal extension
\[
 \pi_0^\#:
 \underline{\mathcal H}^0_{V,W}
 \longrightarrow
 \underline{\operatorname{Hom}}^0(V,W)
\]
therefore gives every sheafified morphism a globally defined smooth
arity-zero component.

\begin{prop}\label{prop:target-whitehead}
Every morphism in $W_{\mathrm{shHom}}$ is invertible in
$H^0(\Loc_{\mathrm{shHom}})$.
\end{prop}

\begin{proof}
Let $f:V\to W$ lie in $W_{\mathrm{shHom}}$.  Its global arity-zero
component $f_0$ is a chain map because the arity-zero part of
$D f=0$ is
\[
 d_Wf_0-f_0d_V=0.
\]
The bounded finite-rank bundle complex $\operatorname{Cone}(f_0)$ is
exact on every fiber.  Starting in its lowest degree, exactness
determines the rank of each differential inductively; all ranks are
locally constant, so kernels and images are smooth subbundles.  Bundle
metrics give global smooth complements, and the differentials restrict
to bundle isomorphisms between the complements and the subsequent
images.  Their inverses define a smooth contracting homotopy.
Consequently $f_0$ is a smooth chain-homotopy equivalence.

For any target object $X$, composition with $f$ gives sheaf maps
\[
 f_*:\underline{\mathcal H}^\bullet_{X,V}
 \longrightarrow\underline{\mathcal H}^\bullet_{X,W},
 \qquad
 f^*:\underline{\mathcal H}^\bullet_{W,X}
 \longrightarrow\underline{\mathcal H}^\bullet_{V,X}.
\]
Filter by arity.  On the associated graded, only $f_0$ contributes:
higher components of $f$ strictly raise arity.  Hence both associated
graded maps are chain-homotopy equivalences.  The arity filtrations are
finite in each total degree, so $f_*$ and $f^*$ are stalkwise
quasi-isomorphisms.  Their cones are bounded-below complexes of fine
sheaves, and thus their maps on global sections are quasi-isomorphisms.

Taking $X=W$ in $f_*$ lifts $[\id_W]$ to a right inverse
$[g_R]\in H^0\operatorname{Hom}(W,V)$, and taking $X=V$ in $f^*$
lifts $[\id_V]$ to a left inverse $[g_L]$ in the same group.  Then
\[
 [g_L]=[g_L][f][g_R]=[g_R],
\]
so the common class is a two-sided inverse to $[f]$.
\end{proof}

This proposition does not construct a global generalized cone for a
sheafified morphism and therefore does not imply that
$\Loc_{\mathrm{shHom}}$ is pretriangulated.

Let $N_{\mathrm{dg}}$ denote the dg nerve and $N_{A_\infty}$ the
$A_\infty$ nerve \cite{Fao17}; for dg nerves and stability, see
\cite[Section 1.3.1]{LurHA}.  Proposition~\ref{prop:coh-localization}
gives
\[
 N_{\mathrm{dg}}(\mathcal P_{\A})
 [W_{\mathrm{coh}}^{-1}]
 \simeq
 N_{\mathrm{dg}}(\mathcal P_{\A}).
\]
Proposition~\ref{prop:target-whitehead} likewise gives
\[
 N_{\mathrm{dg}}(\Loc_{\mathrm{shHom}})
 [W_{\mathrm{shHom}}^{-1}]
 \simeq
 N_{\mathrm{dg}}(\Loc_{\mathrm{shHom}}).
\]
These are localization statements inside the representation
categories; they do not identify their objects or prove essential
surjectivity.

\subsection{The induced functor of infinity-categories}\label{sec:nerve}

A strictly unital $A_\infty$ functor induces a functor of
$A_\infty$ nerves.  Theorem~\ref{thm:integration} therefore gives
\[
 \mathbf{RH}:
 N_{A_\infty}(\mathcal P_{\A})
 \longrightarrow
 N_{A_\infty}\left(
 \Loc_{\mathrm{shHom},\Ch_{\R}}^{\mathrm{sm,perf}}
 (\Pi_\infty(\CF))\right).
\]
Its unary component at a vertex is $f_x^0$, so it preserves the weak
equivalences above and also induces
\[
 N_{A_\infty}(\mathcal P_{\A})[W_{\mathrm{coh}}^{-1}]
 \longrightarrow
 N_{A_\infty}(\Loc_{\mathrm{shHom}})
 [W_{\mathrm{shHom}}^{-1}].
\]

\begin{cor}\label{cor:nerve-rh}
On the cohomologically regular subcategories,
$\mathbf{RH}$ is fully faithful.  If $M$ is compact, it is an
equivalence of infinity-categories.
\end{cor}

\begin{proof}
An $A_\infty$ quasi-fully faithful functor induces equivalences on all
mapping spaces of the $A_\infty$ nerve, and an $A_\infty$
quasi-equivalence induces a categorical equivalence
\cite{Fao17}.  Apply Proposition~\ref{prop:qff} and
Theorem~\ref{thm:compact-rh}.
\end{proof}

\begin{cor}[Hypersheaf form of the correspondence]
\label{cor:hypersheaf-rh}
On the site of open subsets of a regular foliated manifold, set
\begin{align*}
 \mathcal C_{\mathrm{src}}(U)
 &=
 N_{A_\infty}\left(
 \mathcal P_{\Omega^\bullet_{\CF}(U)}^{\mathrm{reg}}
 \right),\\
 \mathcal C_{\mathrm{raw}}(U)
 &=
 N_{A_\infty}\left(
 \Loc_{\mathrm{raw},\Ch_{\R}}^{\mathrm{sm,perf,reg}}
 \bigl(\Pi_\infty(\CF|_U)\bigr)
 \right).
\end{align*}
After a fixed strictification of restriction, raw integration induces
an equivalence of $\operatorname{Cat}_\infty$-valued
hypersheafifications
\[
 \mathcal C_{\mathrm{src}}^{\mathrm{hyp}}
 \simeq
 \mathcal C_{\mathrm{raw}}^{\mathrm{hyp}}.
\]
 Regularity is understood as a local replete condition.  The
 hypersheafification is applied to the displayed presheaves as written
 and may add objects whose amplitude is only locally bounded.
 For every fixed finite amplitude interval, the same statement holds
 after restricting both presheaves before hypersheafification.  We do
 not assert that hypersheafification commutes with the filtered union
 over all amplitude intervals.
 \end{cor}

\begin{proof}
Star-shaped foliated boxes form a basis, and
Theorem~\ref{thm:local-rh} says that raw integration is an equivalence
on every member of this basis.  Hypersheafification is localization at
local equivalences, which may be detected on a basis
\cite[Section 6.5.2]{LurHT}.  The induced morphism of
hypersheafifications is therefore an equivalence.  Restriction
preserves regularity, and the local inverse in the proof of
Theorem~\ref{thm:local-rh} preserves every fixed amplitude interval.
\end{proof}

Corollary~\ref{cor:hypersheaf-rh} is a local statement about the
stackified models.  It does not identify their global sections with
the global raw categories.  In particular, it neither proves nor
assumes Conjecture~\ref{conj:raw-hom}.

Replacing either side by a quasi-equivalent dg or $A_\infty$ model does
not change the resulting infinity-category.  The cohesive category is
pretriangulated, so its dg nerve is stable, and its idempotent
completion is the corresponding idempotent-complete stable
infinity-category.  No parallel pretriangulatedness or stability claim
is made for $\Loc_{\mathrm{shHom}}$: generalized cones are presently
available only for globally raw-representable morphisms.  Its dg nerve
and the explicit localization above nevertheless exist without such a
claim.  Moreover, the regular cohesive full subcategory need not be
closed under cones, so Corollary~\ref{cor:nerve-rh} is not asserted to
be an equivalence of stable infinity-categories.

\subsection{Functoriality and variance}

For a foliated map
\[
 f:(M,\CF_M)\longrightarrow(N,\CF_N),
\]
integration of simplices is covariant:
\[
 \Pi_\infty(f):
 \Pi_\infty(\CF_M)\longrightarrow\Pi_\infty(\CF_N).
\]
Representations, however, pull back contravariantly.  Naturality of
the iterated-integral formulas, together with
Proposition~\ref{prop:sheafified-pullback}, gives the commutative
square below after the fixed strictification (and otherwise a
canonically coherent square):
\[
\begin{tikzcd}[column sep=5.5em]
\mathcal P_{\Omega^\bullet_{\CF_N}(N)}
 \arrow[r,"\RH_{\CF_N}"]
 \arrow[d,"f^*"']
&
\Loc_{\mathrm{shHom},\Ch_{\R}}^{\mathrm{sm,perf}}
 \bigl(\Pi_\infty(\CF_N)\bigr)
 \arrow[d,"\Pi_\infty(f)^*"]
\\
\mathcal P_{\Omega^\bullet_{\CF_M}(M)}
 \arrow[r,"\RH_{\CF_M}"']
&
\Loc_{\mathrm{shHom},\Ch_{\R}}^{\mathrm{sm,perf}}
 \bigl(\Pi_\infty(\CF_M)\bigr).
\end{tikzcd}
\]
The same variance is required in the singular setting considered
below.

\section{Singular foliations and
$L_\infty$-algebroids}\label{sec:singular-outlook}

The preceding constructions suggest an extension to singular
foliations presented by $L_\infty$-algebroids.  We formulate the
completed Chevalley--Eilenberg and representation-theoretic structures
that such an extension would require and describe the expected
comparison.  In this generality, a smooth integration functor, an
appropriate smooth representation theory, and the corresponding
descent theorem remain to be constructed.

Several companion manuscripts in preparation pursue parts of this
singular program.  Block and Zeng study characteristic classes for
singular foliations equipped with $L_\infty$-resolutions
\cite{BlockZengSingularClassesPrep}; Block, Wei, and Zeng develop
global Dolbeault $L_\infty$-algebroid resolutions of singular
holomorphic foliations \cite{BlockWeiZengDolbeaultPrep}; and Zeng
develops higher monodromy models and a geometric program for holonomy
stacks \cite{ZengHigherMonodromyPrep,ZengHolonomyStacksPrep}.  These
works are logically separate from the regular-foliation theorem and
are not used to justify the conjectural constructions below.

\subsection{Completed Chevalley--Eilenberg models}

\begin{defn}\label{def:linfty-algebroid}
Let $\mathfrak g$ be a bounded nonpositively graded vector bundle of
finite rank in every degree, and put
\[
 \widehat S_M\bigl(\mathfrak g^\vee[-1]\bigr)
 =
 \prod_{m\geq0}
 \Gamma\left(
 \operatorname{Sym}^m
 \bigl(\mathfrak g^\vee[-1]\bigr)\right).
\]
Equip this algebra with its symmetric-degree filtration and product
locally convex topology.  An \emph{$L_\infty$-algebroid} structure
on $\mathfrak g$ is a continuous degree-one square-zero derivation
\[
 Q:\widehat S_M\bigl(\mathfrak g^\vee[-1]\bigr)
 \longrightarrow
 \widehat S_M\bigl(\mathfrak g^\vee[-1]\bigr)
\]
whose restriction to $C^\infty(M)$ is a first-order differential
operator with values in the part of symmetric degree one.  Its
Chevalley--Eilenberg algebra is the complete topological cdga
\[
 \CE(\mathfrak g)
 =
 \left(\widehat S_M
 \bigl(\mathfrak g^\vee[-1]\bigr),Q\right).
\]
A cohesive representation of $\mathfrak g$ is a bounded finite-rank
graded vector bundle whose completed extension to
$\CE(\mathfrak g)$ carries a continuous flat $\Z$-connection.
\end{defn}

For modules and representations up to homotopy of Lie
$n$-algebroids, see Jotz Lean--Mehta--Papantonis \cite{JMP23}.  The
completed topological convention above is adapted to the analytic
integration problem considered here; its relation to their algebraic
framework requires a separate comparison.

For an ordinary Lie algebroid concentrated in degree zero, this is its
usual CE algebra.  In particular,
$\CE(\CF)=\Omega^\bullet_{\CF}(M)$ for a regular foliation.  The
completed model isolates the algebraic data relevant to the proposed
extension.  Unlike the regular-foliation complex, it need not carry
the finite form filtration used in the perturbation arguments above.

The relevant integration results in the literature include:
\begin{enumerate}
 \item For compact manifolds, Block--Smith prove the higher
 Riemann--Hilbert quasi-equivalence \cite[Theorem 4.1]{BS14}; their
 $A_\infty$ functor is \cite[Theorem 4.2]{BS14}, quasi-full
 faithfulness is \cite[Proposition 4.3]{BS14}, and essential
 surjectivity is \cite[Theorem 4.5]{BS14}.  Arias
 Abad--Sch\"atz construct the higher $A_\infty$ components and extend
 the construction to representations up to homotopy of every ordinary
 Lie algebroid \cite[Definition 4.16 and Theorem 4.18]{AS12}.
 \item Getzler and Henriques integrate suitable nilpotent or
 finite-type $L_\infty$-algebras to Kan simplicial objects
 \cite{Get09,Hen08}; Rogers--Zhu study the corresponding homotopy
 theory \cite{RZ20}.
 \item A singular foliation admitting a finite geometric resolution
 has a universal $L_\infty$-algebroid, unique up to the appropriate
 homotopy \cite{Lav18,LLS20}.
\end{enumerate}
An extension of the present theorem to an arbitrary
$L_\infty$-algebroid over a manifold requires, in addition, a global
smooth Kan integration, a smooth perfect representation category, and
an $A_\infty$ Riemann--Hilbert functor between them.  Completed
symmetric powers also destroy the finite form-degree filtration that
made the perturbation arguments above automatic.

\subsection{Perfect singular foliations and the expected comparison}

Let $\CF$ now be a singular foliation in the sense of
Androulidakis--Skandalis \cite{AS06}, namely a locally finitely
generated involutive $C^\infty_M$-submodule of vector fields.  Suppose
it admits a finite geometric resolution and choose a universal
$L_\infty$-algebroid $\mathfrak g$.  Write
$\mathcal G_\bullet(\mathfrak g)$ for a putative smooth integration.
The expected comparison sends a higher algebroid simplex to its base
map:
\[
 \tau:
 \mathcal G_\bullet(\mathfrak g)
 \longrightarrow
 \operatorname{Mon}_\bullet(\CF).
\]
A truncation theorem for $\tau$ would require a definition of the
singular monodromy object, smoothness and independence of the chosen
resolution, a homotopical universal property, and a characterization
of the representations that descend through $\tau$.  The existence of
$\mathfrak g$ alone does not imply these properties.

Subject to these constructions, the expected relationship among the
three representation theories is summarized by the following diagram.
In the lower row, $\Loc_{\mathrm{shHom}}$ denotes the analogue of
Proposition~\ref{prop:sheafified-target} whenever the putative
simplicial object has a restriction theory over $M$.  Dashed arrows
denote conjectural extensions in the displayed generality.  The solid
pullback $\tau^*$ is defined once the map $\tau$ and its
target representation category have been constructed.
\[
\begin{tikzcd}[column sep=4.6em,row sep=3.2em]
\mathfrak g
 \arrow[r,dashed,"\operatorname{Int}"]
 \arrow[d,phantom,
 "\scriptstyle\operatorname{Rep}\ \leadsto" description]
&
\mathcal G_\bullet(\mathfrak g)
 \arrow[r,dashed,"\tau"]
 \arrow[d,phantom,
 "\scriptstyle\operatorname{Rep}\ \leadsto" description]
&
\operatorname{Mon}_\bullet(\CF)
 \arrow[d,phantom,
 "\scriptstyle\operatorname{Rep}\ \leadsto" description]
\\
\operatorname{Coh}\bigl(\CE(\mathfrak g)\bigr)
 \arrow[r,dashed,"\RH_{\mathfrak g}"']
&
\Loc_{\mathrm{shHom}}^{\mathrm{sm,perf}}
 \bigl(\mathcal G_\bullet(\mathfrak g)\bigr)
&
\Loc_{\mathrm{shHom}}^{\mathrm{sm,perf}}
 \bigl(\operatorname{Mon}_\bullet(\CF)\bigr)
 \arrow[l,"\tau^*"']
\end{tikzcd}
\]
The vertical symbols indicate representation assignments rather than
covariant functors from geometric objects to categories.  The cited
integration theorems construct the top arrow in their respective
settings, but not in the displayed generality.  Likewise,
$\RH_{\mathfrak g}$ is established for ordinary Lie algebroids and
remains conjectural for general $L_\infty$-algebroids.  The direction
of $\tau^*$ reflects the contravariance of pullback.  A rightward lower
arrow would instead require a derived Kan extension or a descent
functor with additional hypotheses.

Essential surjectivity of $\tau^*$ is not expected without a descent
condition.  For example, a map from a one-object Lie groupoid to a
point pulls point representations back only to trivial group
representations.  This suggests that the essential image of $\tau^*$
should consist of representations whose action along the homotopy
fibers of $\tau$ is coherently trivial.  Constructing $\tau$, proving
such a descent criterion, and extending the reconstruction theorem to
the completed Chevalley--\allowbreak Eilenberg algebra remain open.
These problems are prerequisites for extending the correspondence to
singular foliations and general $L_\infty$-algebroids.

\begingroup
\raggedright
\printbibliography

@article{AC11,
  author={Arias Abad, Camilo and Crainic, Marius},
  title={Representations up to homotopy of {Lie} algebroids},
  journal={Journal f{\"u}r die reine und angewandte Mathematik},
  volume={663},
  year={2012},
  pages={91--126},
  doi={10.1515/CRELLE.2011.095},
  eprint={0901.0319},
  archivePrefix={arXiv},
  primaryClass={math.DG}
}

@article{AMV21,
  author={Arias Abad, Camilo and Pineda Montoya, Santiago and Quintero V{\'e}lez, Alexander},
  title={Chern--Weil theory for {$\infty$}-local systems},
  journal={Transactions of the American Mathematical Society},
  volume={377},
  number={5},
  year={2024},
  pages={3129--3171},
  doi={10.1090/tran/9068},
  eprint={2105.00461},
  archivePrefix={arXiv},
  primaryClass={math.AT}
}

@article{AQVV19,
  author={Arias Abad, Camilo and Quintero V{\'e}lez, Alexander and V{\'e}lez V{\'a}squez, Sebasti{\'a}n},
  title={An {$A_\infty$} version of the {Poincar{\'e}} lemma},
  journal={Pacific Journal of Mathematics},
  volume={302},
  number={2},
  year={2019},
  pages={385--412},
  doi={10.2140/pjm.2019.302.385}
}

@article{AVV21,
  author={Arias Abad, Camilo and V{\'e}lez V{\'a}squez, Sebasti{\'a}n},
  title={The higher {Riemann--Hilbert} correspondence and principal {$2$}-bundles},
  journal={S{\~a}o Paulo Journal of Mathematical Sciences},
  volume={15},
  number={2},
  year={2021},
  pages={637--660},
  doi={10.1007/s40863-021-00239-y}
}

@article{CHL21,
  author={Chuang, Joseph and Holstein, Julian and Lazarev, Andrey},
  title={Maurer--Cartan moduli and theorems of {Riemann--Hilbert} type},
  journal={Applied Categorical Structures},
  volume={29},
  number={4},
  year={2021},
  pages={685--728},
  doi={10.1007/s10485-021-09631-3}
}

@article{Hol15Morita,
  author={Holstein, Julian V. S.},
  title={Morita cohomology},
  journal={Mathematical Proceedings of the Cambridge Philosophical Society},
  volume={158},
  number={1},
  year={2015},
  pages={1--26},
  doi={10.1017/S0305004114000516},
  eprint={1404.1327},
  archivePrefix={arXiv},
  primaryClass={math.AT}
}

@article{Hol15HLC,
  author={Holstein, Julian V. S.},
  title={Morita cohomology and homotopy locally constant sheaves},
  journal={Mathematical Proceedings of the Cambridge Philosophical Society},
  volume={158},
  number={1},
  year={2015},
  pages={27--35},
  doi={10.1017/S0305004114000528},
  eprint={1404.1326},
  archivePrefix={arXiv},
  primaryClass={math.AT}
}

@article{JMP23,
  author={Jotz Lean, Madeleine and Mehta, Rajan Amit and Papantonis, Theocharis},
  title={Modules and representations up to homotopy of {Lie} {$n$}-algebroids},
  journal={Journal of Homotopy and Related Structures},
  volume={18},
  number={1},
  year={2023},
  pages={23--70},
  doi={10.1007/s40062-022-00322-x},
  eprint={2001.01101},
  archivePrefix={arXiv},
  primaryClass={math.DG}
}

@article{TV23,
  author={To{\"e}n, Bertrand and Vezzosi, Gabriele},
  title={Algebraic foliations and derived geometry: the {Riemann--Hilbert} correspondence},
  journal={Selecta Mathematica},
  volume={29},
  number={1},
  year={2023},
  eid={5},
  doi={10.1007/s00029-022-00808-9}
}

@article{AS06,
  author={Androulidakis, Iakovos and Skandalis, Georges},
  title={The holonomy groupoid of a singular foliation},
  journal={Journal f{\"u}r die reine und angewandte Mathematik},
  volume={626},
  year={2009},
  pages={1--37},
  doi={10.1515/CRELLE.2009.001}
}

@article{AS12,
  author = {Arias Abad, Camilo and Schätz, Florian},
  title = {The {$A_\infty$} de {Rham} theorem and integration of representations up to homotopy},
  journal = {International Mathematics Research Notices},
  year = {2013},
  number = {16},
  pages = {3790--3855},
  doi = {10.1093/imrn/rns166}
}

@incollection{Blo05,
  author={Block, Jonathan},
  title={Duality and equivalence of module categories in noncommutative geometry},
  booktitle={A Celebration of the Mathematical Legacy of {Raoul Bott}},
  series={CRM Proceedings and Lecture Notes},
  volume={50},
  pages={311--339},
  publisher={American Mathematical Society},
  location={Providence, RI},
  year={2010},
  doi={10.1090/crmp/050/24},
  eprint={math/0509284},
  archivePrefix={arXiv},
  primaryClass={math.QA}
}

@misc{Blo06II,
  author={Block, Jonathan},
  title={Duality and equivalence of module categories in noncommutative geometry {II}: {Mukai} duality for holomorphic noncommutative tori},
  year={2006},
  eprint={math/0604296},
  archivePrefix={arXiv},
  primaryClass={math.QA},
  url={https://arxiv.org/abs/math/0604296}
}

@article{BenBasBlock13,
  author={Ben-Bassat, Oren and Block, Jonathan},
  title={Milnor descent for cohesive dg-categories},
  journal={Journal of K-Theory},
  volume={12},
  number={3},
  year={2013},
  pages={433--459},
  doi={10.1017/is013007003jkt236},
  eprint={1201.6118},
  archivePrefix={arXiv},
  primaryClass={math.AG}
}

@article{BHW17,
  author={Block, Jonathan and Holstein, Julian V. S. and Wei, Zhaoting},
  title={Explicit homotopy limits of dg-categories and twisted complexes},
  journal={Homology, Homotopy and Applications},
  volume={19},
  number={2},
  year={2017},
  pages={343--371},
  doi={10.4310/HHA.2017.v19.n2.a17},
  eprint={1511.08659},
  archivePrefix={arXiv},
  primaryClass={math.AT}
}

@misc{Cra04,
  author={Crainic, Marius},
  title={On the perturbation lemma, and deformations},
  year={2004},
  eprint={math/0403266},
  archivePrefix={arXiv},
  primaryClass={math.AT},
  doi={10.48550/arXiv.math/0403266}
}

@article{BS14,
  title={The higher {Riemann--Hilbert} correspondence},
  author={Jonathan Block and Aaron M. Smith},
  journal={Advances in Mathematics},
  year={2014},
  volume={252},
  pages={382--405},
  doi={10.1016/j.aim.2013.11.001}
}

@phdthesis{ZengThesis,
  author = {Zeng, Qingyun},
  title = {Derived {Lie} Infinity-Groupoids and Algebroids in Higher
    Differential Geometry},
  institution = {University of Pennsylvania},
  year = {2021},
  eprint = {2210.05856},
  archiveprefix = {arXiv},
  primaryclass = {math.DG},
  url = {https://repository.upenn.edu/entities/publication/02977bc1-814e-42a7-914a-6d0924b1eb8d}
}

@unpublished{BlockZengSingularClassesPrep,
author = {Jonathan Block and Qingyun Zeng},
title = {Characteristic Classes of {$L_\infty$}-Resolved Singular Foliations},
note = {In preparation},
}

@unpublished{ZengTransverseHolomorphicPrep,
author = {Qingyun Zeng},
title = {Derived Categories of Transversely Holomorphic Foliations},
note = {In preparation},
}

@unpublished{BlockZengPenroseWardPrep,
author = {Jonathan Block and Qingyun Zeng},
title = {A Categorical Penrose--Ward Correspondence for Self-Dual Yang--Mills Theory},
note = {In preparation},
}

@unpublished{BlockWeiZengDolbeaultPrep,
author = {Jonathan Block and Zhaoting Wei and Qingyun Zeng},
title = {Global Dolbeault {$L_\infty$} Algebroid Resolutions of Singular Holomorphic Foliations},
note = {In preparation},
}

@unpublished{ZengHigherMonodromyPrep,
author = {Qingyun Zeng},
title = {Higher Monodromy Models for Geometrically Resolved Singular Foliations},
note = {In preparation},
}

@unpublished{ZengHolonomyStacksPrep,
author = {Qingyun Zeng},
title = {Holonomy Stacks of Singular Foliations: A Geometric Program for Higher Transverse Transport},
note = {In preparation},
}

@article{Che77,
  title={Iterated path integrals},
  author={Kuo-Tsai Chen},
  journal={Bulletin of the American Mathematical Society},
  year={1977},
  volume={83},
  number={5},
  pages={831--879},
  doi={10.1090/S0002-9904-1977-14320-6}
}

@article{Gug77,
author = {V. K. A. M. Gugenheim},
title = {{On Chen's iterated integrals}},
volume = {21},
journal = {Illinois Journal of Mathematics},
number = {3},
publisher = {Duke University Press},
pages = {703 -- 715},
year = {1977},
doi = {10.1215/ijm/1256049021}
}

@article{Hen08,
author = {Henriques, André},
year = {2008},
month = {07},
pages = {1017--1045},
number = {4},
title = {Integrating {$L_\infty$}-algebras},
volume = {144},
journal = {Compositio Mathematica},
doi = {10.1112/S0010437X07003405}
}

@misc{Igu09,
	title={Iterated integrals of superconnections},
	author={Kiyoshi Igusa},
	year={2009},
	eprint={0912.0249},
	archivePrefix={arXiv},
	primaryClass={math.AT},
	url={https://arxiv.org/abs/0912.0249}
}

@phdthesis{Lav18,
  author={Lavau, Sylvain},
  title={Lie {$\infty$}-algebroids and singular foliations},
  institution={Universit{\'e} Claude Bernard Lyon 1},
  type={PhD thesis},
  year={2016},
  url={https://theses.hal.science/tel-01447963},
  eprint={1703.07404},
  archivePrefix={arXiv},
  primaryClass={math.DG}
}

@article{LLS20,
  title={The universal {Lie} infinity-algebroid of a singular foliation},
  author={Camille Laurent-Gengoux and Sylvain Lavau and Thomas Strobl},
  journal={Documenta Mathematica},
  year={2020},
  volume={25},
  pages={1571--1652},
  doi={10.4171/DM/782}
}

@book{Gro55,
  author={Grothendieck, Alexandre},
  title={Produits tensoriels topologiques et espaces nucl{\'e}aires},
  series={Memoirs of the American Mathematical Society},
  number={16},
  publisher={American Mathematical Society},
  address={Providence, RI},
  year={1955},
  doi={10.1090/memo/0016}
}

@book{Tre67,
  author={Tr{\`e}ves, Fran{\c c}ois},
  title={Topological Vector Spaces, Distributions and Kernels},
  series={Pure and Applied Mathematics},
  volume={25},
  publisher={Academic Press},
  address={New York},
  year={1967}
}

@article{Wei23,
  author={Wei, Zhaoting},
  title={Descent of dg cohesive modules for open covers on complex manifolds},
  journal={European Journal of Mathematics},
  volume={9},
  number={3},
  year={2023},
  eid={77},
  doi={10.1007/s40879-023-00675-4},
  eprint={1804.00993},
  archivePrefix={arXiv},
  primaryClass={math.AG}
}

@article{Fao17,
  author={Faonte, Giovanni},
  title={Simplicial nerve of an {$A_\infty$}-category},
  journal={Theory and Applications of Categories},
  volume={32},
  number={2},
  year={2017},
  pages={31--52},
  doi={10.70930/tac/868236x7}
}

@misc{LurHA,
  author={Lurie, Jacob},
  title={Higher Algebra},
  year={2017},
  url={https://www.math.ias.edu/~lurie/papers/HA.pdf},
  note={Online manuscript}
}

@book{LurHT,
  author={Lurie, Jacob},
  title={Higher Topos Theory},
  series={Annals of Mathematics Studies},
  volume={170},
  publisher={Princeton University Press},
  location={Princeton, NJ},
  year={2009},
  isbn={9780691140490},
  doi={10.1515/9781400830558}
}

@book{May67,
  author={May, J. Peter},
  title={Simplicial Objects in Algebraic Topology},
  series={Van Nostrand Mathematical Studies},
  volume={11},
  publisher={D. Van Nostrand Company},
  location={Princeton, NJ},
  year={1967}
}

@article{Get09,
  author={Getzler, Ezra},
  title={Lie theory for nilpotent {$L_\infty$}-algebras},
  journal={Annals of Mathematics},
  volume={170},
  number={1},
  year={2009},
  pages={271--301},
  doi={10.4007/annals.2009.170.271}
}

@article{RZ20,
  author={Rogers, Christopher L. and Zhu, Chenchang},
  title={On the homotopy theory for Lie infinity-groupoids, with an application to integrating {$L_\infty$}-algebras},
  journal={Algebraic \& Geometric Topology},
  volume={20},
  number={3},
  year={2020},
  pages={1127--1219},
  doi={10.2140/agt.2020.20.1127}
}
\endgroup

\end{document}